\documentclass[a4paper]{amsart}
\usepackage{amscd}

\newcommand{\bs}{\backslash}
\newcommand{\longto}{\longrightarrow}
\newcommand{\R}{\mathbb{R}}
\newcommand{\C}{\mathbb{C}}
\newcommand{\w}{\omega}
\newcommand{\s}{\sigma}
\newcommand{\pphi}{\varphi}
\newcommand{\taum}{\tau^{-}}
\newcommand{\qhamsp}{(M,\w,\mu:M\to U)}
\newcommand{\qhamquot}{\mu^{-1}(\{1\})/U}
\newcommand{\normK}{\mathcal{N}(K)}
\newcommand{\Walc}{\mathcal{W}}
\newcommand{\clWalc}{\overline{\mathcal{W}}}
\newcommand{\im}{\mathrm{Im}}

\newcommand{\Mgl}{\mathcal{M}_{g,l}}
\newcommand{\Mod}{\mathrm{Hom}_{\mathcal{C}}(\pi_{g,l},U)/U}
\newcommand{\Hom}{\mathrm{Hom}}
\newcommand{\Rep}{\mathrm{Rep}}
\newcommand{\calC}{\mathcal{C}}
\newcommand{\bhat}{\hat{\beta}}
\newcommand{\Mext}{(U\times U)\times\cdots\times (U\times U)\times\mathcal{C}_1\times\cdots\times \mathcal{C}_l}
\newcommand{\Mtot}{(U\times U)^g \times\mathcal{C}_1\times\cdots\times \mathcal{C}_l}
\newcommand{\Sig}{\Sigma}
\newcommand{\pigl}{\pi_1(\Sigma_g\backslash\{s_1,\, ...\, ,s_l\})}
\newcommand{\prespiintro}{<\alpha_1,\beta_1,\, ...\, ,\alpha_g,\beta_g,\gamma_1,\, ...\, ,\gamma_l\ |\ \prod_{i=1}^{g} [\alpha_i,\beta_i] \prod_{j=1}^{l} \gamma_j =1>}
\newcommand{\prespi}{(a_1,b_1,\, ...\, ,a_g,b_g,c_1,\, ...\, ,c_l)\in U^{2g+l}\ |\ \prod_{i=1}^g [a_i,b_i] \prod_{j=1}^l c_j =1}
\newcommand{\abc}{(a_1,b_1,\, ...\, ,a_g,b_g,c_1,\, ...\, ,c_l)}
\newcommand{\relabc}{[a_1,b_1]...[a_g,b_g]c_1...c_l}
\newcommand{\relabcshort}{\prod_{i=1}^g [a_i,b_i] \prod_{j=1}^l c_j}
\newcommand{\pconj}{\mathcal{C}_1\times\cdots\times\mathcal{C}_l}
\newcommand{\HomC}{\mathrm{Hom}_{\mathcal{C}}}
\newcommand{\RepC}{\mathrm{Rep}_{\mathcal{C}}}
\newcommand{\piS}{\pi_1(S^2\backslash\{s_1,\, ...\,  ,s_l\})}
\newcommand{\tU}{\widetilde{U}}
\newcommand{\tu}{\widetilde{u}}
\newcommand{\tM}{\widetilde{M}}
\newcommand{\tw}{\widetilde{\omega}}
\newcommand{\tmu}{\widetilde{\mu}}
\newcommand{\tbeta}{\widetilde{\beta}}
\newcommand{\ttau}{\widetilde{\tau}}
\newcommand{\ttaum}{\widetilde{\tau}^{-}}
\newcommand{\tQ}{\widetilde{Q}}
\newcommand{\Z}{\mathcal{Z}}
\newcommand{\tX}{\widetilde{X}}

\title[Decomposable representations of surface groups]{Decomposable representations and Lagrangian submanifolds of moduli spaces associated to surface groups}

\thanks{Supported by the Japanese Society for Promotion of Science (JSPS)}

\newtheorem{theorem}{Theorem}[section]
\newtheorem{lemma}[theorem]{Lemma}
\newtheorem{proposition}[theorem]{Proposition}
\newtheorem{corollary}[theorem]{Corollary}
\newtheorem{definition}[theorem]{Definition}
\newtheorem{remark}[theorem]{Remark}
\newtheorem*{thm*}{Theorem}
\newtheorem*{ack}{Acknowledgments}

\author{Florent Schaffhauser}

\address{Keio University\\ 
Dept. of Mathematics\\
Hiyoshi 3-14-1\\
Kohoku-ku, 223-8522\\
Yokohama, Japon}

\email{florent@math.keio.ac.jp}



\begin{document}

\begin{abstract}
The importance of explicit examples of Lagrangian submanifolds of moduli spaces is revealed by papers such as \cite{DS,Salamon}:  given a $3$-manifold $M$ with boundary $\partial M=\Sigma$, Dostoglou and Salamon use such examples to obtain a proof of the Atiyah-Floer conjecture relating the symplectic Floer homology of the representation space $\Hom(\pi_1(\Sigma=\partial M), U)/U$ (associated to an explicit pair of Lagrangian submanifolds of this representation space) and the instanton homology of the $3$-manifold $M$. In the present paper, we construct a Lagrangian submanifold of the space of representations $\Mgl:=\Mod$  of the fundamental group $\pi_{g,l}$ of a punctured Riemann surface $\Sigma_{g,l}$ into an arbitrary compact connected Lie group $U$. This Lagrangian submanifold is obtained as the fixed-point set of an anti-symplectic involution $\bhat$ defined on $\Mgl$. We show that the involution $\bhat$ is induced by a form-reversing involution $\beta$ defined on the quasi-Hamiltonian space $\Mtot$. The fact that $\bhat$ has a non-empty fixed-point set is a consequence of the real convexity theorem for group-valued momentum maps proved in \cite{Sch_IM}. The notion of decomposable representation provides a geometric interpretation of the Lagrangian submanifold thus obtained.
\end{abstract}

\maketitle

\section{Introduction}

The purpose of this article is to give examples of Lagrangian submanifolds in the moduli spaces 
$$\Mgl:=\Mod$$ 
where $$\pi_{g,l}:=\prespiintro$$
is the fundamental group of a Riemann surface $\Sig_{g,l}:=\Sig_g\bs\{s_1,\, ...\, ,s_l\}$ with $g$ holes ($g\geq 0$ is the genus of the compact Riemann surface $\Sig_g$) and $l$ punctures ($l\geq 0$), and $U$ is an arbitrary compact connected Lie group. Let us explain how we intend to proceed.\\
The above-mentioned moduli spaces $\Mgl=\Mod$  are the spaces of (equivalence classes of) representations of
the fundamental group $\pi_{g,l}:=\pi_1(\Sigma_{g,l})$ of a Riemann surface
$\Sigma_{g,l}:=\Sigma_g\bs\{s_1,\, ... \, ,s_l\}$ where $\Sigma_g$ is
a compact Riemann surface of genus $g\geq 0$, where $l$ is an integer
$l\geq 0$ (with the convention that $\Sigma_{g,0}:=\Sigma_g$) and where
$s_1,\, ...\, ,s_l$ are $l$ pairwise distinct points of $\Sigma_g$.
These \emph{representations varieties} have been an important object of study
for several decades now, and are located at the intersection of various
areas of mathematics, each of which is very rich and sheds interesting
light on these spaces. Thus, the space
$$\Rep(\pi_{g,l},U):=\Hom(\pi_{g,l},U)/U$$ of equivalence classes of
representations of $\pi_{g,l}$ in a Lie group $U$ arises naturally in
complex algebraic geometry as it can be identified to the space of
equivalence classes of holomorphic vector bundles on $\Sigma_{g,l}$, as was
shown by Narasimhan and Seshadri in the 1960s (see \cite{Nara-Se}). At the
beginning of the 1980s,
Atiyah and Bott gave a new direction to the subject (see \cite{AB}) by 
identifying these spaces as the moduli spaces of flat connections on
principal bundles of group $U$ on $\Sigma_{g,l}$, thereby revealing the
importance of the representation varieties in gauge theory.\\
In this paper, we chose to investigate the \emph{symplectic} structure of these representation spaces. This symplectic structure can be obtained and described in a variety of ways
(see for instance \cite{Gol,GHJW,Al-Mal1,AMM,MW}), each of which has its own
advantages. The description given by Alekseev, Malkin
and Meinrenken in \cite{AMM} will prove particularly well-suited for
our study of representations of $\pi_{g,l}$. This description rests on the
notion of \emph{quasi-Hamiltonian space}, which enables one to avoid
infinite-dimensional manifolds while limiting oneself to relatively simple
objects to construct a symplectic form on representation varieties.\\ 
To be able to be more precise in our statements, let us recall that the group $\pi_{g,l}=\pigl$ admits the following finite presentation by generators and relations:
$$\pi_{g,l}=\prespiintro.$$
In this presentation, $\alpha_1,\beta_1,\, ...\, ,\alpha_g,\beta_g$ stand for homotopy classes of non-trivial loops in the compact Riemann surface $\Sig_g$, and $\gamma_1,\, ...\, ,\gamma_l$ stand for homotopy classes of loops around the punctures $s_1,\, ...\, ,s_l$ ($\gamma_j$ wraps one time only around $s_j$ and not around any other $s_k$). Let $U$ be a compact connected Lie group. We denote by 
$$\Hom(\pi_{g,l},U):=\{\rho:\pi_{g,l}\longto U\}$$
the set of group morphisms of $\pi$ into $U$. Elements of $\Hom(\pi_{g,l},U)$ are also called \emph{representations of} $\pi_{g,l}$ \emph{into} $U$. By choosing generators of $\pi_{g,l}$, one may identify $\Hom(\pi_{g,l},U)$ with a subset of $U^{2g+l}$:
$$\Hom(\pi_{g,l},U)$$ $$\simeq
\{\prespi\}\subset U^{2g+l}.$$
Two representations $\rho=\abc$ and $\rho'=(a'_1,b'_1,\, ...\, ,$ $a'_g,b'_g,c'_1,\,$  $...\, c'_l)$ of $\pi_{g,l}$ into $U$ are said to be \emph{equivalent} if there exists an element $\pphi\in U$ such that $a'_i=\pphi a_i \pphi^{-1}, b'_i=\pphi b_i\pphi^{-1}$ and $c'_j=\pphi c_j\pphi^{-1}$ for all $i$ and $j$. In other words, the representations $\rho$ and $\rho'$ are equivalent if they are in a same orbit of the diagonal conjugation action of $U$ on $U^{2g+l}$ (since this action preserves the relation $\prod[a_i,b_i]\prod c_j=1$, it leaves $\Hom(\pi_{g,l},U)$ stable). When $U\subset Gl(V)$ is a group of linear transformations, this is indeed the same notion as equivalence of linear representations. In order to obtain symplectic structures on the orbit space $\Hom(\pi_{g,l},U)/U$ of this action (that is, the space of equivalence classes of representations), one has to fix the conjugacy classes of the generators $c_1,\, ...\, ,c_l$ \emph{representing homotopy classes of loops around the punctures} $s_1,\, ... ,s_l$ (see for instance \cite{AMM}). Otherwise, one only obtains a Poisson structure (see for instance \cite{AKSM}). Observe that the conjugacy classes of the $a_i,b_i$ need \emph{not} be fixed. Let us now denote by $\calC_1,\, ...\, ,\calC_l$ a collection of $l$ conjugacy classes of $U$ and study the set 
$$\Hom_\calC(\pi_{g,l},U)$$ $$:=
\{\abc \in (U\times U)^g \times\calC_1\times\cdots\times\calC_l\ | \prod_{i=1}^{g} [a_i,b_i] \prod_{j=1}^l c_j=1\}$$
$$\subset (U\times U)^g \times \pconj$$
which is a (possibly empty) subset of the set $\Hom(\pi_{g,l},U)\subset (U\times U)^g\times U^l$. In all of the following, we will assume that the conjugacy classes $\calC_1,\, ...\, ,\calC_l$ of $U$ are chosen in a way that $\HomC(\pi_{g,l},U)\not=\emptyset$. When $g=0$ (that is, for the case of the punctured sphere group), giving necessary and sufficient conditions for this to be true is a difficult problem (see for instance \cite{AW}). When $g\geq 1$, however, and if $U$ is semi-simple, one always has $\HomC(\pi_{g,l},U)\not=\emptyset$ (see \cite{Ho-Liu2}). Since each conjugacy class $\calC_j$ is preserved by the conjugacy action, the set $\HomC(\pi_{g,l},U)$ is stable under the diagonal action of $U$. Then, the \emph{representation variety} (or \emph{representation space}) is by definition the orbit space 
$$\RepC(\pi_{g,l},U):=\Hom_C(\pi_{g,l},U)/U$$
of this action. From now on, we will call \emph{representations of} $\pi_{g,l}$ \emph{into} $U$ the elements of $\HomC(\pi_{g,l},U)$ (meaning that we will assume that $c_1,\, ...\, ,c_l$ lie in prescribed conjugacy classes). To avoid confusion, we will denote the space of equivalence classes of such representations by 
$$\Mgl:=\Mod$$
and call it the \emph{moduli space associated to} $\pi_{g,l}$. For a generic choice of conjugacy classes $(\calC_j)_{1\leq j\leq l}$, the moduli spaces $\Mgl$ ($l\geq 1$) are smooth symplectic manifolds. Here, by a generic choice we mean conjugacy classes of maximal dimension: $\calC_j$ is the conjugacy class of an element $c_j\in U$ whose centralizer $\mathcal{Z}(c_j)$ is a maximal torus of $U$ (if $U=U(n)$, this means that $c_j$ has pairwise distinct eigenvalues). For special choices of conjugacy classes however, these moduli spaces are \emph{stratified symplectic spaces} in the sense of Lerman and Sjamaar (see \cite{LS} and \cite{GHJW}). To us, this will simply mean that $\Mgl$ is a disjoint union of
smooth manifolds (of different dimensions) called strata, each of which carries a symplectic
structure. Recall from \cite{AMM} that the moduli space $\Mgl$ is the quasi-Hamiltonian quotient
$$\Mgl=\mu^{-1}(\{1\})/U$$
associated to the quasi-Hamiltonian space $M:=(U\times U)^g\times \pconj$ with momentum map 
\begin{eqnarray*}
\mu: (U\times U)\times\cdots\times(U\times U)\times\pconj & \longto & U\\
\abc & \longmapsto & \prod_{i=1}^g [a_i,b_i] \prod_{j=1}^lc_j
\end{eqnarray*}
which is one possible way of seeing that $\Mgl$ carries a symplectic structure (we shall come back to this in section \ref{involutions}). A quasi-Hamiltonian quotient is often denoted by:
$$M//U:=\mu^{-1}(\{1\})/U.$$
As mentioned earlier, our objective is to give an example of a Lagrangian submanifold in the moduli space $\Mgl$ for any $g,l\geq 0$ and for an arbitrary compact connected Lie group $U$. While the meaning of the term \emph{Lagrangian submanifold} is clear when $\Mgl$ is a smooth symplectic manifold, the possibly singular analogue still lacks, as far as we know, a precise definition when $\Mgl$ is a stratified symplectic space. Let us propose the following one:
\begin{definition}[Stratified Lagrangian subspace of a stratified symplectic space]\label{def_strat_lag_sub}
Let $N=\sqcup_{j\in J} X_j$ be a disjoint union of symplectic manifolds (for instance, $N=M//U$ a quasi-Hamiltonian quotient). A subset $L\subset N$ satisfying:\begin{enumerate}
\item[(i)] $L\not=\emptyset$.
\item[(ii)] if $L\cap X_j\not=\emptyset$ then $L\cap X_j$ is a Lagrangian submanifold of the symplectic manifold $X_j$. 
\end{enumerate}
is called a \emph{stratified Lagrangian subspace of} $N$.
\end{definition}
Thus, a \emph{stratified Lagrangian subspace of} $\Mgl$ is a non-empty subset of
$\Mgl$ whose intersection with each stratum, if non-empty, is a Lagrangian submanifold of
the considered stratum. Since the moduli spaces $\Mgl$ ($l\geq 1$) are smooth symplectic manifolds for a generic choice of conjugacy classes $(\calC_j)_{1\leq j\leq l}$, we will often, when speaking informally, talk about a Lagrangian submanifold of $\Mgl$. In the statements of our results, however, we will consistently identify which objects are smooth symplectic manifolds and which objects are stratified symplectic spaces. In the latter case, we will use the terminology \emph{stratified Lagrangian subspace} to refer to the possibly singular analogue of a a Lagrangian submanifold in that context (see definition \ref{def_strat_lag_sub}). The path we will follow to find a Lagrangian submanifold of the moduli space $\Mgl=\Mod$ consists in:
\begin{enumerate}
\item[1.] introducing a notion of decomposable  representation of $\pi_{g,l}$ into $U$.
\item[2.] characterizing these representations in terms of an involution $\beta$ defined on $M=\Mtot$.
\item[3.] showing that this involution $\beta$ induces an involution $\bhat$ on the moduli space $\Mgl=\Mod=M//U$ whose fixed-point is a Lagrangian submanifold of $\Mgl$.
\end{enumerate}
In section \ref{decomp}, we will define decomposable representations of $\pi_{g,l}$ into $U$ (definition \ref{decomp_rep}). Since the notion of decomposable representations is somewhat unfamiliar, let us mention here that it is distinct from the notion of \emph{reducible} or \emph{irreducible} representation: a given representation of $\pi_{g,l}$ into $U$ might be decomposable or not regardless of whether it is irreducible. As a matter of fact, an interesting question is to ask whether there always exist, for given $\pi_{g,l}$ and $U$, irreducible representations of $\pi_{g,l}$ into $U$ which are decomposable. As we shall see in subsection \ref{decomp_and_irred}, this turns out to be always true when ($g=0$ and $l\geq 1$) or ($g\geq 1$ and $l=0$), and generically true when ($g\geq 1$ and $l\geq1$). Having defined decomposable representations,
our task will then be twofold: first, we will have to obtain a characterization of these decomposable representations in terms of an involution $\beta$ defined on the quasi-Hamiltonian space $M=\Mtot$ (theorem \ref{charac_thm}) and then we will have to prove that decomposable representations always exist (theorem \ref{existence_thm_III}). The existence problem is the more difficult of the two: the proof that decomposable representations always exist rests on their characterization in terms of an involution and on a real convexity theorem for group-valued momentum maps obtained in \cite{Sch_IM}. More precisely, recall that the moduli space $\Mgl$ is the quasi-Hamiltonian quotient $$\Mgl=M//U=\mu^{-1}(\{1\})/U$$ where $M$ is the quasi-Hamiltonian space $M=(U\times U)^g \times\pconj$, with momentum map $$\mu\abc=\prod_{i=1}^g [a_i,b_i]\prod_{j=1}^l c_j.$$ We will define an involution $\beta$ on $M$ such that a representation $$(a,b,c)\in \mu^{-1}(\{1\})\subset (U\times U)^g\times\pconj$$ is decomposable if and only if $\beta(a,b,c)=\pphi.(a,b,c)$ for some \emph{symmetric} element $\pphi$ of $U$ (with respect to a fixed involutive automorphism $\tau$ of $U$, see theorem \ref{charac_thm}), and we will show that $\beta$ satisfies sufficient conditions for it to induce an involution $\bhat$ on $\Mgl=\mu^{-1}(\{1\})/U$ (see proposition \ref{lag_locus}). We will then use the real convexity theorem of \cite{Sch_IM} to prove that $\mu^{-1}(\{1\})\cap Fix(\beta)\not=\emptyset$, which will imply that $\bhat$ has a non-empty fixed-point set. Combining this with the fact that $\beta$ reverses the $2$-form defining the quasi-Hamiltonian structure on $M$ (proposition \ref{beta_anti-inv}) we will finally prove that $Fix(\bhat)$ is a stratified Lagrangian subspace of the moduli space $\Mgl=\mu^{-1}(\{1\})/U$. The main result of this paper is theorem \ref{lag_nature_thm}, which we now state to summarize the above:
\begin{theorem}[A stratified Lagrangian subspace of $\Mgl$]
For any compact connected Lie group $U$ and any integers $g,l\geq 0$, the involution $\beta: M \to M$ introduced in definition \ref{def_beta} induces an involution $\bhat$ on the moduli space $\Mgl=\Mod$ of representations of the surface group $\pi_{g,l}=\pigl$ into $U$ whose fixed-point set $Fix(\bhat)$ is a stratified Lagrangian subspace of $\Mgl$. In particular, for a generic choice of conjugacy classes $(\calC_j)_{1\leq j\leq l}$, the set $Fix(\bhat)$ is a Lagrangian submanifold of the smooth symplectic manifold $\Mgl$ ($l\geq 1$).
\end{theorem}
In section \ref{application}, the last one of this paper, we will give an application of the results obtained when characterizing decomposable representations of the fundamental group $\piS$ of a punctured sphere in terms of an involution $\beta$ (defined, in this case, on the quasi-Hamiltonian space $\pconj$).   The result we will obtain there (theorem \ref{compact_Thompson}) can be thought of as an analogue of Thompson's problem in the case of a compact connected Lie group. We will underline the analogy with the symplectic geometry of Thompson's problem as it is treated in \cite{AMW}.\\
Other examples of Lagrangian submanifolds of representation spaces are given in \cite{Goldman_moduli} and \cite{Ho}. In this last reference \cite{Ho}, Ho also obtains such a Lagrangian submanifold by constructing an anti-symplectic involution on the moduli space $\mathcal{M}_{g,0}=\HomC(\pi_{g,0},U)/U$ of representations of the fundamental group of a \emph{compact} Riemann surface of genus $g\geq 1$.

\section{The notion of decomposable representation}\label{decomp}

In this section, we give the definition of decomposable representations of surface groups (see definition \ref{decomp_rep}). We begin with a special case: the case of the fundamental group $\pi_{0,l}$ of an $l$-punctured sphere (case $g=0$). The purpose of doing so is to stress that the notion of decomposable representations when $g\geq 1$ is derived from the case where $g=0$. In the latter case, this notion has a simple geometric origin and was first introduced by Falbel and Wentworth in their study of representations of the fundamental group $\pi_{g,0}$ of a punctured sphere into the unitary group $U(n)$ (see \cite{FW}). Their approach (for the group $\pi_{g,0}$) was generalized to arbitrary compact connected Lie groups $(U,\tau)$ endowed with an involutive automorphism in \cite{Sch_CJM}. After reviewing this in subsection \ref{punctured_sphere}, we will introduce the notion of decomposable representation for surface groups $\pi_{g,l}$ with $g\geq 1$ in subsection \ref{higher_genus}. 

\subsection{The case of the punctured sphere group}\label{punctured_sphere}

As mentioned above, the notion of decomposable representation in the case of the fundamental group of the punctured sphere is due to Falbel and Wentworth when $U=U(n)$ (see \cite{FW}, where decomposable representations are called \emph{Lagrangian representations}). For this notion to make sense in an arbitrary compact connected Lie group $U$, we assume that this group $U$ is endowed with an involutive automorphism $\tau$. In fact, we shall make a further assumption on $\tau$, which will be crucial to obtain the characterization of decomposable representations in terms of an involution $\beta$ (see remark \ref{max_rk_assump_I}) as well as to prove their existence (see remark \ref{max_rk_assump_II}). Namely, we will assume that the involution $\tau$ satisfies the following condition:
\begin{definition}[Involutive automorphism of maximal rank]\label{max_rk_inv}
Let $U$ be a compact connected Lie group. An involutive automorphism $\tau$ of $U$ is said to be \emph{of maximal rank} if there exists a \emph{maximal} torus $T\subset U$ fixed \emph{pointwise} by the involution $\taum(u):=\tau(u^{-1})$.
\end{definition}
\begin{proposition}
Let $U$ be a compact connected Lie group. Then there always exists an involutive automorphism of maximal rank of $U$.
\end{proposition}
\begin{proof}
We refer to \cite{Loos} for the proof: the involution $\tau$ is the restriction to $U$ of the involution on the complexification $U^{\C}$ of $U$ defining the non-compact dual $G=Fix(\tau)\subset U^{\C}$ of $U$. The term maximal rank refers to the fact that the symmetric space $U/Fix(\tau|_U)$ is of maximal rank.
\end{proof}
\begin{remark}\label{isometry}
In the following, the Lie algebra $\mathfrak{u}=Lie(U)=T_1U$ of our compact connected Lie group $U$ will always be endowed with an $Ad$-invariant scalar product $(.\, |\, .)$. We will assume that this scalar product is chosen such that $T_1\tau:\mathfrak{u}\to\mathfrak{u}$ is an isometry.
\end{remark}
As an example, one can think of $(U=U(n),\tau(u)=\overline{u})$, so that $\taum(u)=u^t$ and the maximal torus of $U(n)$ consisting of diagonal unitary matrices is fixed pointwise by $\taum$. In this case, $U^{\C}=Gl(n,\C)$ and $G=Gl(n,\R)$. As a $Ad$-invariant scalar product on $\mathfrak{u}(n)$, one may choose $(X|Y)=-tr(XY)$, for which complex conjugation is an isometry. Observe that if $\tau$ is an arbitrary involutive automorphism of a compact connected Lie group $U$, there always exists a torus fixed pointwise by $\taum$ but it is not necessarily a \emph{maximal} torus (hence the term maximal rank, in the latter case). In the following, we will only consider \emph{symmetric pairs} $(U,\tau)$ such that $\taum$ satisfies the assumptions of definition \ref{max_rk_inv} above. In particular, we have:
\begin{lemma}\label{restriction}
If $\tau$ is an involutive automorphism of maximal rank of $U$, then the involution $\taum: u\in U\mapsto \tau(u^{-1})$ sends a conjugacy class $\calC\subset U$ into itself.
\end{lemma}
\begin{proof}
This follows from the fact that there exists a maximal torus $T$ of $U$ which is
pointwise fixed by $\taum$. Indeed, if we take
$u\in\calC$, then $u=vdv^{-1}$ for some $v\in U$ and some $d\in T$. Then
$\taum(u)=\tau(v)\taum(d)\taum(v)=\tau(v)d\taum(v)$ is conjugate to $d$ and
therefore to $u$.
\end{proof}
\begin{definition}[Decomposable representations of $\piS$]\label{zero-genus_decomp_rep}\ 
Let $(U,\tau)$ be a compact connected Lie group endowed with an involutive automorphism $\tau$ of maximal rank. A representation $(c_1, \, ...\, c_l)$ of the group
$\pi_{0,l}=\piS$ into $U$ is called \emph{decomposable} if there exist $l$ elements
$w_1,\, ...\, ,w_l$ $\in U$ satisfying:
\begin{enumerate}
\item[(i)] $\tau(w_j)=w_j^{-1}$ for all $j$ (each $w_j$ is a \emph{symmetric}
element of $U$ with respect to $\tau$).
\item[(ii)] $c_1= w_1 w_2^{-1}$, $c_2=w_2 w_3^{-1}$,\, ...\, , $c_l=w_l w_1^{-1}$.
\end{enumerate}
\end{definition}
For details on the geometric origin of this definition, we refer to \cite{Sch_CJM}: when $U=U(n)$, the decomposition $u_j=w_jw_{j+1}^{-1}$ is equivalent to saying that $u_j$ is a product $u_j=\sigma_j\sigma_{j+1}$ of orthogonal reflections with respect to Lagrangian subspaces of $\C^n$ (hence the name Lagrangian representation in \cite{FW}). Recall that $$\pi_{0,l}=\piS=<\gamma_1,\, ...\, ,\gamma_l\ |\ \gamma_1...\gamma_l=1>$$ and observe that if $c_1,\, ...\, ,c_l$ satisfy condition (ii) in definition \ref{zero-genus_decomp_rep} above, they automatically satisfy $c_1...c_l=1$. Condition (i) requires \emph{additionally} that each generator $c_j$ be a product of two \emph{symmetric} elements: $c_j=w_jw_{j+1}^{-1}$ with $\taum(w_j)=w_j$ and $\taum(w_{j+1})=w_{j+1}$.

\subsection{The case of an arbitrary surface group}\label{higher_genus}

We will now derive the definition of decomposable representation for an arbitrary surface group $\pi_{g,l}$ from definition \ref{zero-genus_decomp_rep} above. Recall that \begin{eqnarray*}\pi_{g,l} & = & \pigl\\
 & = & \prespiintro.\end{eqnarray*}
The key to obtaining the notion of a decomposable representation of $\pi_{g,l}$ is to observe that if $\abc$ is a representation of the group $\pi_{g,l}=\pigl$ then $([a_1,b_1],\, ...\, ,[a_g,b_g],c_1, $ $...\, ,c_l)$ is a representation of the fundamental group $\pi_{0,g+l}=\pi_1(S^2\bs\{r_1,\, ...\, ,r_g,s_1,\, ...\, ,s_l\})$ of a $(g+l)$-punctured sphere: $$\underbrace{[a_1,b_1]}_{=u_1}...\underbrace{[a_g,b_g]}_{=u_g}\underbrace{c_1}_{=u_{g+1}}...\underbrace{c_l}_{=u_{g+l}}=1.$$ Building on this remark, we will begin by asking, in addition to $c_j=w_jw_{j+1}^{-1}$ with $w_j,w_{j+1}\in Fix(\taum)$, that each $[a_i,b_i]$ be decomposable under the form $[a_i,b_i]=v_iv_{i+1}^{-1}$ with $v_i,v_{i+1}\in Fix(\taum)$ for all $i$. For consistency reasons, we will also ask that $v_{g+1}=w_1$ and $w_{l+1}=v_1$. As one might suspect, this is not restrictive enough; we need additional conditions on the $(a_i,b_i)$. Namely, we will ask:
\begin{definition}[Decomposable representations of $\pigl$]\label{decomp_rep}
Let $(U,\tau)$ be a compact connected Lie group endowed with an involutive automorphism $\tau$ of maximal rank. A representation $(a_1,b_1,\, ...\, ,a_g,b_g,c_1,\, ...\, ,$ $c_l)$ of $\pigl$ into $U$ is called \emph{decomposable} if there exist $(g+l)$ elements $v_1,\, ...\, ,v_g,w_1, \, ...\, , w_l \in U$ satisfying:
\begin{enumerate}
\item[(i)] $\tau(v_i)=v_i^{-1}$ for all $i$ and $\tau(w_j)=w_j^{-1}$ for all $j$. 
\item[(ii)] $[a_1,b_1]=v_1 v_2^{-1}$, $[a_2,b_2]=v_2 v_3^{-1}$,\, ...\, , $[a_g,b_g]=v_g w_1^{-1}$, $c_1= w_1 w_2^{-1}$, $c_2=w_2 w_3^{-1}$,\, ...\, , $c_l=w_l v_1^{-1}$.
\item[(iii)] $\tau(a_i)=v_{i+1}^{-1}b_iv_{i+1}$ for all $i\in \{1,\, ...\, ,g\}$ (with $v_{g+1}=w_1$).
\end{enumerate}
\end{definition}
To understand where condition (iii) comes from, let us take a look at the case $l=0$. Then in particular we have, from (ii), that $[a_g,b_g]=v_g v_1^{-1}$. To simplify even further, let  us consider the special case $v_1=1\in U$. Then one has $[a_g,b_g]\in Fix(\taum)$, that is: $a_gb_ga_g^{-1}b_g^{-1}=\tau(b_g)\tau(a_g)\tau(b_g^{-1})\tau(a_g^{-1})$. One may then observe that this is always true if $\tau(a_g)=b_g$. Condition (iii) may then be derived from the above analysis by lifting the simplifying assumptions made. As a matter of fact, \emph{in the} $g\geq 1$ \emph{case}, the above definition was obtained also by searching for an involution $\beta$ on $\Mtot$ and not only by looking for an appropriate notion of decomposable representation. This might give the impression that the above notion of decomposable representation is an \emph{ad hoc} notion. To some extent, that is indeed the case. As we shall see later, the key object in this paper is the involution $\beta$ on $\Mtot$ defined in definition \ref{def_beta}, since it is this involution that will induce an anti-symplectic involution $\bhat$ on $\Mgl$, whose fixed-point set $Fix(\bhat)$ will provide the promised Lagrangian submanifold. And this involution $\beta$ was obtained by trying to generalize the expression found in the $g=0$ case, for which we refer to \cite{Sch_CJM}. But since in this $g=0$ case, the involution $\beta$ was obtained \emph{from} the notion of decomposable representation, it seemed natural to work out a notion of decomposable representation in the $g\geq 1$ case, that will later provide a nice geometric interpretation of the Lagrangian submanifold $Fix(\bhat)\subset \Mgl$. One may notice that definition \ref{decomp_rep} coincides with definition \ref{zero-genus_decomp_rep} when $g=0$.
\begin{remark}
Observe that in definition \ref{decomp_rep} above, we do not ask each of the generators $a_i$ and $b_i$ (which account for the surface having non-zero genus, and whose conjugacy class is \emph{not} prescribed) to be decomposable as a product of two symmetric elements.
\end{remark}
I would like to thank Pierre Will for numerous discussions ultimately leading to definition \ref{decomp_rep}. For another notion of decomposable representation, we refer to \cite{Will}.
Before ending this section, we observe the following useful fact that the notion of decomposable representation is compatible with the equivalence of representations:
 \begin{proposition}\label{equiv_decomp_rep}
 Let $\abc$ be a representation of $\pi_{g,l}=\pigl$ and $u$ be any element of $U$. Then the representation $$\abc$$ is decomposable if and only if the representation $$u.\abc:=(u a_1u^{-1}, u b_1u^{-1},\, ...\, , u c_lu^{-1})$$ is decomposable.
 \end{proposition}
 \begin{proof}
 If $\abc$ is decomposable, then there exist $$v_1,\, ...\, ,v_g,w_1,\, ...\, ,w_l \in Fix(\taum)$$ satisfying $[a_i,b_i]=v_i v_{i+1}^{-1}$, $c_j=w_j w_{j+1}^{-1}$ and $\tau(a_i)=v_{i+1}^{-1}b_iv_{i+1}$ for all $i,j$ (with $v_{g+1}=w_1$ and $w_{l+1}=v_1$). Let us note $a'_i:=u a_iu^{-1}$, $b'_i=u b_iu^{-1}$ and $c'_j=u c_ju^{-1}$. Then:
 \begin{eqnarray*}
 [a'_i,b'_i] & = & u [a_i,b_i] u^{-1}\\
 & = & u v_iv_{i+1}^{-1} u^{-1}\\
 & = & (u v_i \taum(u))(\tau(u)v_{i+1}^{-1}u^{-1})\\
 & = & (u v_i \taum(u))(u v_{i+1}\taum(u))^{-1}\\
 & = & v'_i (v'_{i+1})^{-1}
\end{eqnarray*}
with $v'_i=u v_i \taum(u) \in Fix(\taum)$, and likewise:
 \begin{eqnarray*}
 c'_j & = & u c_j u^{-1}\\
 & = & (u w_j \taum(u))(u w_{j+1}\taum(u))^{-1}\\
 & = & w'_j (w'_{j+1})^{-1}
 \end{eqnarray*}
with $w'_j=u w_j \taum(u)\in Fix(\taum)$. Finally:
\begin{eqnarray*}
\tau(a'_i) & = & \tau(u) \tau(a_i)\taum(u)\\
& = & \tau(u) v_{i+1}^{-1} b_i  v_{i+1} \taum(u)\\
& = & \tau(u) v_{i+1}^{-1} (u^{-1}u) b_i (u^{-1}u) v_{i+1} \taum(u)\\
& = & (v'_{i+1})^{-1} b'_i v'_{i+1}
\end{eqnarray*}
which proves that the representation $(a'_1,b'_1,\, ...\, ,c'_l)$ is decomposable.\\
The converse implication is proved the same way, setting $v_i:=u^{-1}v'_i\tau(u)$ and $w'_j:=u^{-1}w'_j\tau(u)$.
\end{proof}
Our task will now be, as announced in the introduction, to obtain a characterization of decomposable representations of $\pi_{g,l}$ in terms of an involution $\beta$ defined on the quasi-Hamiltonian space $\Mtot$. This will be achieved in theorem \ref{charac_thm}. But first, let us give the generalities about anti-symplectic involutions on quasi-Hamiltonian quotients that we will need in the following.
 
\section{Anti-symplectic involutions on quasi-Hamiltonian quotients}\label{involutions}

In this section, we will recall the generalities on quasi-Hamiltonian spaces that we shall need in the following. The reference for the notion of quasi-Hamiltonian space is the original paper of Alekseev, Malkin and Meinrenken (\cite{AMM}). For further details on the structure of a quasi-Hamiltonian quotient $M//U=\mu^{-1}(\{1\})/U$, we refer to \cite{Sch_Pekin}. Some of the results we are about to state also come from \cite{Sch_CJM}.

\subsection{Structure of quasi-Hamiltonian quotients}

Let us begin by recalling the definition of a quasi-Hamiltonian space, which is due to Alekseev, Malkin and Meinrenken in \cite{AMM}.
\begin{definition}[Quasi-Hamiltonian space, \cite{AMM}]
Let $(M,\w)$ be a manifold endowed with a $2$-form $\w$ and an action of the
Lie group $(U,(.\,|\, .))$ leaving the $2$-form $\w$ invariant. We denote by
$(.\,|\, .)$ an $Ad$-invariant non-degenerate symmetric bilinear form on
$\mathfrak{u}:=Lie(U)$, by $\theta^L:=u^{-1}.du$ and $\theta^R:=du.u^{-1}$ the Maurer-Cartan $1$-forms of $U$, and by $\chi:=\frac{1}{2}([\theta^L,\theta^L]\, |\, \theta^L)$ the Cartan $3$-form of $U$. Finally, we denote by $X^{\sharp}$ the fundamental vector field on $M$ associated to $X\in\mathfrak{u}$. Its value at $x\in M$ is : $X^{\sharp}_x:=\frac{d}{dt}|_{t=0}(\exp(tX).x)$.  Let $\mu:M
\to U$ be a $U$-equivariant map (for the conjugacy action of $U$ on itself).\\ Then
$(M,\w,\mu:M\to U)$ is said to be a \emph{quasi-Hamiltonian space} (with respect to the action of $U$) if the map $\mu: M\to U$ satisfies the following three conditions:
\begin{enumerate}
\item[(i)] $d\w = -\mu^* \chi$.
\item[(ii)] for all $x \in M$, $\ker \w_x= \{X^{\sharp}_x\ :\
X\in\mathfrak{u}\,|\, (Ad\,\mu(x)+Id).X=0\}$.
\item[(iii)] for all $X\in \mathfrak{u}$, $\iota_{X^{\sharp}}\w =
\frac{1}{2}\mu^*(\theta^L +\theta^R \,|\, X)$.
\end{enumerate}
\noindent where $(\theta^L+\theta^R\,|\, X)$ is the real-valued $1$-form
defined on $U$ for any $X\in\mathfrak{u}$ by $(\theta^L + \theta ^R \,|\,
X)_u(\xi) := (\theta^L_u(\xi) + \theta^R_u(\xi)\,|\, X)$ (where $u\in U$ and
$\xi \in T_u U$). The map $\mu$ is called the \emph{momentum map}. 
\end{definition}
Let us now study, given a quasi-Hamiltonian $U$-space $(M,\w,\mu:M\to U)$, the associated quasi-Hamiltonian quotient $M//U:=\mu^{-1}(\{1\}) /U$, which may also be denoted by $M^{red}$. By a theorem of Alekseev, Malkin and Meinrenken, this space admits a symplectic structure. More precisely, it is shown in \cite{AMM} that if $1\in U$ is a regular value of the momentum map $\mu$, then $M//U$ is a symplectic orbifold. Using the same techniques and combining with the notion of stratified symplectic space, one obtains the following result:
\newcommand{\calN}{\mathcal{N}}
\begin{proposition}[Structure of a quasi-Hamiltonian quotient, \cite{Sch_Pekin}]\label{quotientstructure}
Let $(M,\w,\mu:M\to U)$ be a quasi-Hamiltonian $U$-space. For any closed subgroup $K\subset U$, denote by $M_K$ the isotropy manifold of type $K$ in $M$: $$M_K=\{x\in M\ |\ U_x=K\}.$$ Denote by $\calN(K)$ the normalizer of $K$ in $U$ and by $L_K$ the quotient group $L_K:=\calN(K)/K$. Then $\mu(M_K)\subset\calN(K)$ and if we denote by $\mu_K$ the composed map $\mu:M_K\to \calN(K)\to L_K=\calN(K)/K$, then $(M_K,\w|_{M_K},\mu_K:M_L\to L_K)$ is a quasi-Hamiltonian $L_K$-space. Furthermore, $L_K$ acts freely on $M_K$ and the orbit space $$(\mu^{-1}(\{1_U\})\cap M_K)/L_K$$ is a symplectic manifold.\\ Denote by $(K_j)_{j\in J}$ a system of representatives of closed subgroups of $U$. Then the orbit space $M^{red}:=\mu^{-1}(\{1_U\})/U$ is the disjoint union of the following symplectic manifolds:
$$\mu^{-1}(\{1_U\})/U = \bigsqcup_{j\in J} (\mu^{-1}(\{1_U\})\cap M_{K_j})/L_{K_j}$$
\end{proposition}\newcommand{\fiber}{\mu^{-1}(\{1\})}
In particular, if $U$ acts freely on $\mu^{-1}(\{1\})$, the quasi-Hamiltonian quotient $M//U=\fiber/U$ is a symplectic manifold, as shown in \cite{AMM}. In the following, a subset $X_K:=(\fiber\cap M_K)/L_K$ in the above description of $M//U$ will be called a \emph{stratum}. We end this subsection by recalling the following result:
\begin{lemma}\label{classical_lemma}
Let $(N,\w)$ be a symplectic manifold and $\s$ be an anti-symplectic
involution on $N$ (meaning that $\s^*\w=-\w$ and $\s^2 = Id_N$).
Denote by $N^{\s}:=Fix(\s)$ the fixed-point set of $\s$. Then: if
$N^{\s}\not=\emptyset$, it is a Lagrangian submanifold of $N$.
\end{lemma}
Observe that an anti-symplectic involution does not necessarily have fixed
points. For instance, the map $(-Id_{\R^3})|_{S^2}: (x,y,z)\in S^2 \mapsto
-(x,y,z)$ reverses orientation on $S^2$ (so that it is anti-symplectic with
respect to the
volume form $x\, dy\wedge dz -y\, dx\wedge dz -z\, dx\wedge dy$ on $S^2$), and has
no fixed points on $S^2$.

\subsection{Involutions induced on quasi-Hamiltonian quotients}\label{anti-symp_inv}

We now would like to apply lemma \ref{classical_lemma} to the symplectic manifold $N=M//U$, where $\qhamsp$ is a quasi-Hamiltonian $U$-space. The idea is to obtain the involution $\s$ on $N=M//U$ from an involution $\beta$ defined on the quasi-Hamiltonian space $M$. Drawing from the usual Hamiltonian case studied in \cite{OSS}, we set:
\begin{definition}[Compatible involutions]\label{compa}
Let $(M,\w,\mu:M\to U)$ be a quasi-Hamiltonian space and let $\tau$ be an
involutive automorphism of $U$. Denote by $\taum$ the involution on $U$
defined by $\taum(u)=\tau(u^{-1})$. An involution $\beta$ on $M$ is said to
be \emph{compatible with the action of} $U$ if $\beta(u.x)=\tau(u).\beta(x)$
for all $x\in M$ and all $u\in U$, and it is said to be \emph{compatible
with the momentum map} $\mu$ if $\mu\circ\beta=\taum\circ\mu$.
\end{definition}
We then have the following result:
\begin{proposition}[\cite{Sch_CJM}]\label{lag_locus}
Let $(M,\w,\mu:M\to U)$ be a quasi-Hamiltonian space and let $\tau$ be an
involutive automorphism of $U$. Denote by $\taum$ the involution on $U$
defined by $\taum(u)=\tau(u^{-1})$ and let $\beta$ be an involution on
$M$ compatible with the action of $U$ and with the momentum map $\mu$.
Then $\beta$ induces an involution $\hat{\beta}$ on the
quasi-Hamiltonian quotient $M//U:=\mu^{-1}(\{1\})/U$ defined by $\bhat([x])=[\beta(x)]$. If in addition $U$ acts freely on $\fiber$ and $\beta^* \w = -\w$, then $M//U$ is a symplectic manifold and $\hat{\beta}$ is an anti-symplectic involution. Consequently, if $\bhat$ has fixed points, then $Fix(\hat{\beta})$ is a Lagrangian submanifold of $M//U$.
\end{proposition}
\begin{remark}[On the assumption that $Fix(\bhat)\not=\emptyset$]\label{fix_bhat}
We will see in section \ref{existence} that the assumption $Fix(\bhat)\not=\emptyset$ in proposition \ref{lag_locus} above is in fact \emph{always} satisfied, for $U$ a compact connected and \emph{semi-simple} Lie group, provided the involution $\beta$ on $M$ has fixed points whose image lies in the connected component of $Fix(\taum)\subset U$ containing $1$ (see theorem \ref{inv_semi-simple}). For the quasi-Hamiltonian space $M=(U\times U)^g\times\pconj$ with $U$ an \emph{arbitrary} compact connected Lie group, the involution $\beta$ we will consider (see definition \ref{def_beta}) also induces an involution $\bhat$ that always has fixed points (see theorem \ref{existence_thm_III}).
\end{remark}
We will now drop the assumption that $U$ acts freely on $\fiber$ and give more details on the structure of $Fix(\bhat)$ using the description of the quasi-Hamiltonian quotient $M//U$ as the disjoint union of the symplectic manifolds $X_K=(\fiber \cap M_K)/L_K$ over  closed subgroups $K$ of $U$ (see proposition \ref{quotientstructure}). More precisely, we will show that $Fix(\bhat)$ is a union of Lagrangian submanifolds of some of the strata. The first question we ask ourselves is: under what conditions does a stratum $X_K=(\fiber\cap M_K)/L_K\subset M//U$ contain points of $Fix(\bhat)$? We then observe the following facts:
\begin{lemma}\label{stab_beta_x}
Let $\beta$ be a form-reversing involution on the quasi-Hamiltonian space $\qhamsp$ compatible with $\tau:U\to U$ and $\mu:M\to U$. Then for all $x\in M$ the stabilizer $U_{\beta(x)}$ of $\beta(x)$ in $U$ is the subgroup $U_{\beta(x)}=\tau(U_x)$.
\end{lemma}
\begin{proof}
If $x\in U_{\beta(x)}$ then $\beta(\tau(u).x)=u.\beta(x)=\beta(x)$, hence $\tau(u).x=x$. Therefore $\tau(u)\in U_x$ and $u=\tau(\tau(u))\in \tau(U_x)$. Conversely, if $u\in U_x$, then $\beta(x)=\beta(u.x)=\tau(u).\beta(x)$ hence $\tau(u)\in U_{\beta(x)}$.
\end{proof}
\begin{lemma}\label{conjugate_subgroups}
If $[x]\in (\fiber\cap M_K)/L_K$ satisfies $\bhat([x])=[x]$ then the subgroups $K$ and $\tau(K)$ are conjugate in $U$. Moreover, $\bhat$ maps a stratum $X_K$ to itself if and only if $\tau(K)$ is conjugate to $K$ in $U$ (a condition we will denote $\tau(K)\sim K$).
\end{lemma}
\begin{proof}
If $x\in M_K$ satisfies $\bhat([x])=[x]$, we have $U_x=K$ and $\beta(x)=u.x$ for some $u\in U$. Hence $U_{\beta(x)}=uU_xu^{-1}=uKu^{-1}$. By lemma \ref{stab_beta_x}, we also have $U_{\beta(x)}=\tau(U_x)=\tau(K)$, hence $\tau(K)=uKu^{-1}$. Finally, $\bhat$ maps $X_K$ to itself if and only if for all $y\in \fiber\cap M_K$ the stabilizer of $\beta(y)$ is conjugate to $K$, that is: $U_{\beta(y)}=\tau(U_y)=\tau(K)$ is conjugate to $K$.
\end{proof}
\begin{lemma}\label{K_and_MK_stable}
Let  $\beta:M\to M$ be an involution on $M$ compatible with $\tau:U\to U$ and let $K\subset U$ be a (closed) subgroup of $U$. Then $\beta(M_K)\subset M_K$ if and only if $\tau(K)\subset K$. 
\end{lemma}
\begin{proof}
This is a direct consequence of the fact that $U_{\beta(x)}=\tau(U_x)$.
\end{proof}
Observe that since $\beta$ and $\tau$ are involutions, the previous conditions are in fact equivalent to $\beta(M_K)=M_K$ and $\tau(K)=K$. The following statement is then immediate:
\begin{proposition}\label{fixbetadescription}
Let $$M//U=\sqcup_{j\in J} (\fiber\cap M_{K_j})/L_{K_j}$$ be a quasi-Hamiltonian quotient. Set $X_j=(\fiber\cap M_{K_j})/L_{K_j}$. Then $\bhat(X_j)=X_j$ if and only if $\tau(K_j)$ is conjugate to $K_j$ in $U$ and one has: $$Fix(\bhat)=\bigsqcup_{j\in J~|~\tau(K_j)\sim K_j} Fix(\bhat|_{X_j}).$$ In particular, $Fix(\bhat)$ is a disjoint union of Lagrangian submanifolds of the strata $X_j\subset \fiber/U$ for which $\tau(K_j)$ is conjugate to $K_j$ and $Fix(\bhat|_{X_j})\not=\emptyset$.
\end{proposition}
Observe that it is not clear which strata actually satisfy $Fix(\bhat|_{X_j})\not=\emptyset$ although there is at least one such $j$ since $Fix(\bhat)$ is asssumed to be non-empty (see also proposition \ref{ppal_stratum}). The second question we ask ourselves is: does the involution $\bhat|_{X_j}$ come from a form-reversing involution on the quasi-Hamiltonian space $M_{K_j}$? The next result shows that the answer is yes if $\tau(K_j)$ is conjugate to $K_j$ \emph{by an element} $w_j\in Fix(\taum)$ (that is, by a \emph{symmetric} element of $U$ with respect to $\tau$). We will use this remark in subsection \ref{interpretation} (proposition \ref{almost_there}).
\begin{proposition}\label{symmetric_case}
Let $K$ be a closed subgroup of $U$ such that $\tau(K)=w_KKw_K^{-1}$ with $w_K\in Fix(\taum)\subset U$. Then, there exists an involution $\beta_K:M\to M$ and an involutive automorphism $\tau_K:U\to U$ such that $\beta_K$ is compatible with $\tau_K$ and with the momentum map $\mu:M\to U$, and such that $\beta_K(M_K)=M_K$. Moreover, the involution $\widehat{\beta_K}$ induced by $\beta_K$ on $\mu^{-1}(\{1\})/U$ coincides with $\bhat$. In particular on $X_K=(\fiber\cap M_{K})/L_K$, one has: $Fix(\bhat|_{X_K})=Fix(\widehat{\beta_{K}})$. Finally, if $\beta$ is form-reversing then so is $\beta_K$.
\end{proposition}
\begin{proof}
Set $\beta_K(x):=w_K^{-1}.\beta(x)$ for all $x\in M$ and $\tau_K(u)=(w_K)^{-1}\tau(u)w_K$ for all $u\in U$. Using the fact that $w_K$ is a symmetric element of $U$, one easily checks that $\beta_K$ is an involution on $M$, compatible with $\tau_K$ and $\mu$, which induces the same involution as $\beta$ on $\fiber/U$. Further, we have $\tau_K(K)=(w_K)^{-1}\tau(K)w_K=K$. Therefore, by lemma \ref{K_and_MK_stable}, we have $\beta_K(M_K)=M_K$. Finally, $\beta_K=w_K^{-1}.\beta$ reverses the $2$-form $\w$ because $\beta$ reverses $\w$ and $w_K$ leaves it invariant.
\end{proof}
We then observe:
\begin{lemma}\label{NK_stable}
If $\tau_0$ is any involutive automorphism of $U$ and if $K\subset U$ is a subgroup such that $\tau_0(K)\subset K$, one also has $\tau_0(\normK)\subset \normK$. Consequently, $\tau_0$ induces an involution on $L_K=\normK/K$.
\end{lemma}
\begin{proof}
Consider $n\in\normK$ and $k\in K$. Since $\tau_0(k)\in K$, we have $n\tau_0(k)n^{-1}\in K$, and therefore also $\tau_0(n)k(\tau_0(n))^{-1}=\tau_0(n\tau_0(k)n^{-1})\in \tau_0(K)\subset K$, thus $\tau_0(n)\in\normK$.
\end{proof}
It is then immediate that the involution $(\beta_K)|_{M_K}:M_K\to M_K$ is compatible with the action of $(L_K,\tau_K)$ and with the momentum map of this action. For additional details on the structure of the fixed-point set of an anti-symplectic involution on a Hamiltonian or quasi-Hamiltonian quotient, we refer to \cite{Foth_loci} and \cite{Sch_Poisson2006}.

\subsection{Constructing form-reversing involutions on product spaces}

We end this section with a result that will be useful in sections \ref{charac} and forward. It gives a way of constructing form-reversing involutions on product of quasi-Hamiltonian spaces starting from form-reversing involutions on each factor.
\begin{proposition}[\cite{Sch_CJM}]\label{anti-inv}
Let $(M_1,\w_1,\mu_1: M_1 \to U)$ and $(M_2,$ $\w_2,\mu_2: M_2 \to U)$ be two
quasi-Hamiltonian $U$-spaces. Let $\tau$ be an involutive automorphism of
$(U,(. \,|\, .))$ and let $\beta_i$ be an involution on $M_i$ satisfying:
\begin{enumerate}
\item[(i)] $\beta_i^* \w_i = -\w_i$.
\item[(ii)] $\beta_i(u.x_i) = \tau(u).\beta_i(x_i)$ for all $u \in U$ and all
$x_i \in M_i$.
\item[(iii)] $\mu_i \circ \beta_i = \taum \circ \mu_i$.
\end{enumerate}
\noindent Consider the quasi-Hamiltonian $U$-space $$(M:= M_1 \times M_2, \w:= \w_1
\oplus \w_2 + \frac{1}{2}(\mu_1^*\theta^L \wedge \mu_2^* \theta^R), \mu:=\mu_1 \cdot
\mu_2)$$ (with respect to the diagonal action of $U$) and the map:
\begin{eqnarray*}
\beta:= \big( (\mu_2 \circ \beta_2).\beta_1, \beta_2 \big): M & \longto & M \\
(x_1,x_2) & \longmapsto & \big( (\mu_2\circ\beta_2(x_2)).\beta_1(x_1),
\beta_2(x_2)\big)
\end{eqnarray*}
\noindent Then $\beta$ is an involution on $M$ satisfying:
\begin{enumerate}
\item[(i)] $\beta^* \w = -\w$.
\item[(ii)] $\beta(u.x) = \tau(u).\beta(x)$ for all $u \in U$ and all
$x \in M$.
\item[(iii)] $\mu \circ \beta = \taum \circ \mu$.
\end{enumerate}
Moreover, if $\beta_1$ and $\beta_2$ have fixed points, then so does $\beta$.
\end{proposition}
\begin{remark}\label{for_later}
It is even true that if $$x_1\in (\mu_1^{-1}(\{1\})\cap Fix(\beta_1))$$ and $$x_2\in (\mu_2^{-1}(\{1\})\cap Fix(\beta_2))$$ then $$(x_1,x_2)\in (\mu^{-1}(\{1\})\cap Fix(\beta))$$ but in general points of $(\fiber\cap Fix(\beta))$ are not of this form: for instance, if $M=\calC_1\times\calC_2$ a product of two non-trivial conjugacy classes, one has $\mu_i^{-1}(\{1\})=\emptyset$ but this does not necessarily mean that $\mu^{-1}(\{1\})=\emptyset$.
\end{remark}
We now move on to obtaining a characterization of decomposable representations of $\pi_{g,l}$ in terms of an involution $\beta$ satisfying the assumptions encountered in the present section.

\section{Characterization of decomposable representations}\label{charac}

In this section, we give a characterization of decomposable representations of $$\pi_{g,l}=\pigl$$ into a compact connected Lie group $U$ in terms of an involution defined on the quasi-Hamiltonian space $M=\Mtot$ (see \cite{AMM}) whose momentum map
\begin{eqnarray*}
\mu: (U\times U)\times\cdots\times(U\times U)\times\pconj & \longto & U\\
\abc & \longmapsto & \prod_{i=1}^g [a_i,b_i] \prod_{j=1}^lc_j
\end{eqnarray*}
satisfies $$\HomC (\pi_{g,l},U)/U=\qhamquot.$$
As in section \ref{decomp}, the Lie group $U$ is endowed with an involutive automorphism $\tau$ of maximal rank (see definition \ref{max_rk_inv}). The idea of characterizing decomposable representations in terms of an involution comes from the study of the special case $U=U(n)$, as explained in detail in \cite{Sch_CJM}, a case where one can develop a more geometric intuition of the problem using Lagrangian subspaces of $\C^n$ instead of symmetric elements of $U(n)$. The case of the punctured sphere was treated in \cite{Sch_CJM}. We now obtain a similar result for arbitrary surface groups.

\subsection{The involution $\beta$}

We begin by introducing the following map $\beta$:
\begin{definition}[The involution $\beta$]\label{def_beta}
Let $(U,\tau)$ be a compact connected Lie group endowed with an involutive automorphism of maximal rank and let $\calC_1,\, ...\, ,\calC_l$ be $l$ conjugacy classes of $U$. We then define the following map $\beta$ from $\Mtot$ to itself:
\begin{eqnarray*}
\beta: \Mtot & \longto & \Mtot \\
(a_1,b_1,\, ...\, ,a_g,b_g,c_1,\, ...\, ,c_l) & \longmapsto & \Big(\taum([a_2,b_2]...c_l)\tau(b_1)\tau([a_2,b_2]...c_l),\\
& &\quad \taum([a_2,b_2]...c_l)\tau(a_1)\tau([a_2,b_2]...c_l),\\
& & \qquad \qquad \qquad \vdots\\
& & \qquad\taum(c_1...c_l)\tau(b_g)\tau(c_1...c_l),\\
& & \qquad\taum(c_1...c_l)\tau(a_g)\tau(c_1...c_l),\\
& & \qquad\taum(c_2...c_l)\taum(c_1)\tau(c_2...c_l),\\
& & \qquad \qquad \qquad\vdots\\
& & \qquad\qquad\qquad \taum(c_l)\Big)
\end{eqnarray*}
\end{definition}
\begin{remark}
For detailed explanation on how this map was obtained in the $g=0$ case, we refer to \cite{Sch_CJM}. The above generalization to an arbitrary $g\geq 0$ comes from the expression obtained there and the search for an involution $\beta$ on $\Mtot$ satisfying the compatibility conditions of definition \ref{compa}. We will show in this section (theorem \ref{charac_thm}) that the involution $\beta$ of definition \ref{def_beta} admits a very nice geometric interpretation, namely that it characterizes representations of $\pigl$ which are decomposable in the sense of definition \ref{decomp_rep}.
\end{remark}
The above expression for $\beta$ is quite complicated, so we will explain it a little. The term $\beta(a_1,b_1,\, ...\,$ $ a_g,b_g,c_1\, ...\, ,c_l)$ is an element of $\Mtot$. Its first coordinate is conjugate to $\tau(b_1)\in U$ via the element $\taum([a_2,b_2]...[a_g,b_g]c_1...c_l)$, its second coordinate is conjugate to $\tau(a_1)$ via that same element, and so on until its $(2g)^{th}$ coordinate, which is conjugate to $\tau(a_g)\in U$ via $\taum(c_1...c_l)$. Its $(2g+1)^{th}$ coordinate is conjugate to $\taum(c_1)\in\calC_1$ via the element $\taum(c_2...c_l)$, and so on until its $(2g+l)^{th}$ coordinate, which is equal to $\taum(c_l)$. This calls for two remarks.
\begin{remark}
The fact that $\beta$ is well-defined is a consequence of the fact that $\tau$ is of maximal rank: the involution $\taum$ sends a conjugacy class of $U$ into itself (see lemma \ref{restriction}).
\end{remark}
\begin{remark}
Observe that the $a_i,b_i$, which account for the surface having non-zero genus, switch place when they are applied the map $\beta$. This is not a problem because the conjugacy classes of the $a_i,b_i$ are \emph{not} prescribed.
\end{remark}
To make things clearer yet, we now specialize $\beta$ to some particular surface groups.
\begin{itemize}
\item $g=0,l=1$: in this case, $\beta: c\in\calC\mapsto\taum(c)$. We then observe:
\begin{lemma}\label{anti-inv_on_conj_class}
The involution $\beta=\taum|_{\calC}:c\in \calC\mapsto \tau(c^{-1})$ reverses the $2$-form $\w_{\calC}$ defining the quasi-Hamiltonian structure on $\calC$, that is: $\beta^*\w_{\calC}=-\w_{\calC}$.
\end{lemma}
\begin{proof} We refer to \cite{Sch_CJM}.
\end{proof}
\item $g=0,l=3$: (\cite{Sch_CJM}) $$\beta(c_1,c_2,c_3)=(\taum(c_2c_3)\taum(c_1)\tau(c_2c_3),\taum(c_3)\taum(c_2)\tau(c_3),\taum(c_3)).$$
\item $g=1,l=0$: in this case, $\beta: (a,b)\in U\times U \mapsto(\tau(b),\tau(a))$. We then observe:
\begin{lemma}\label{anti-inv_on_double}
The involution $\beta:(a,b)\in U\times U\mapsto (\tau(b),\tau(a))$ reverses the $2$-form $\w_D$ defining the quasi-Hamiltonian structure on $U\times U$, that is: $\beta^*\w_D=-\w_D$.
\end{lemma}
\begin{proof} The $2$-form $\w_D$ defining the quasi-Hamiltonian structure on the double $U\times U$ is $$\w_D = \frac{1}{2}(a^*\theta^L\wedge b^*\theta^R) + \frac{1}{2}(a^*\theta^R\wedge b^*\theta^L) + \frac{1}{2}((a.b)^*\theta^L\wedge (a^{-1}.b^{-1})^*\theta^R)$$ (see \cite{AMM}). Therefore, if we take $(v_i,w_i)_{i=1,2}\in T_aU\times T_bU$, we obtain:
\begin{eqnarray*}
& & (\w_D)_{(a,b)}.\big((v_1,w_1),(v_2,w_2)\big)\\
& = & \frac{1}{2} \big((a^{-1}.v_1\ |\ w_2.b^{-1}) - (a^{-1}.v_2\ |\ w_1.b^{-1})\big)\\
& & + \frac{1}{2} \big((v_1.a^{-1}\ |\ b^{-1}.w_2) - (v_2.a^{-1}\ |\ b^{-1}.w_1)\big)\\
& & + \frac{1}{2} \big(Ad\ a^{-1}.(b^{-1}.w_1) + a^{-1}.v_1\ |\ b^{-1}.w_2 + Ad\ b^{-1}.(a^{-1}.v_2)\big) \\
& & - \frac{1}{2} \Big\{a\leftrightarrow b, v_i\leftrightarrow w_i \Big\}
\end{eqnarray*}
where the notation $a\leftrightarrow b$ means that we switch $a$ and $b$ in the previous line expression (and likewise for $v_i$ and $w_i$). To lighten notation, we denote $\tau(u$) by $\overline{u}$. Then $\beta(a,b)=(\overline{b},\overline{a})$ and:
\begin{eqnarray*}
& & (\beta^*\w_D)_{(a,b)}\big((v_1,w_1),(v_2,w_2)\big) \\
& = & \frac{1}{2} \big( (\overline{b}\, ^{-1}.\overline{w_1}\ |\ \overline{v_2}.\overline{a}\, ^{-1}) - (\overline{b}\, ^{-1}.\overline{w_2}\ |\ \overline{v_1}.\overline{a}\, ^{-1})\big)\\
& & + \frac{1}{2} \big( (\overline{w_1}.\overline{b}\, ^{-1}\ |\ \overline{a}\, ^{-1}.\overline{v_2}.) - (\overline{w_2}.\overline{b}\, ^{-1}\ |\ \overline{a}\, ^{-1}.\overline{v_1})\big)\\
& & - \frac{1}{2} \big( Ad\ \overline{a}\, ^{-1} .(\overline{b}\, ^{-1}.\overline{w_1}) + \overline{a}\, ^{-1}.v_1\ |\ \overline{b}\, ^{-1}.\overline{w_2} + Ad\ \overline{b}\, ^{-1}.(\overline{a}\, ^{-1}.\overline{v_2})\big)\\
& &  + \frac{1}{2} \Big\{a\leftrightarrow b, v_i\leftrightarrow w_i \Big\}
\end{eqnarray*}
Since $\tau$ is an isometry for $(.\ |\ .)$ (see remark \ref{isometry}), we can eliminate it in the above expression for $\beta^*\w_D$ and compare the resulting expression with $\w_D$: we obtain $\beta^*\w_D=-\w_D$.
\end{proof}
\item $g=2,l=0$:
\begin{eqnarray*}
& & \beta(a_1,b_1,a_2,b_2)\\ 
& = & \big(\taum([a_2,b_2])\tau(b_1)\tau([a_2,b_2]), \taum([a_2,b_2])\tau(a_1)\tau([a_2,b_2]), 
\tau(b_2), \tau(a_2) \big)
\end{eqnarray*}
\item $g=2,l=2$:
\begin{eqnarray*}
\beta(a_1,b_1,a_2,b_2) & = & \Big(\taum([a_2,b_2]c_1c_2)\tau(b_1)\tau([a_2,b_2]c_1c_2),\\ 
& & \qquad\qquad\quad ...\, \tau(a_1)\, ... ,\\
& &\qquad \taum(c_1c_2)\tau(b_2)\tau(c_1c_2)\\
& & \qquad\qquad ...\, \tau(a_2)\, ...\, ,\\ 
& & \qquad\taum(c_2)\taum(c_1)\tau(c_2),\\
& & \qquad\qquad\quad \taum(c_2) \Big)
\end{eqnarray*}
\end{itemize}
The fact that $\beta$ is an involution can be checked directly, along with compatibility with the diagonal action of $U$ on $\Mtot$ and the momentum map $\mu\abc=\relabcshort$. But, as the most important thing will be to prove that $\beta$ reverses the $2$-form defining the quasi-Hamiltonian structure on $M=\Mtot$, we will show that we can obtain all of these properties by applying proposition \ref{anti-inv} to the map $\beta$. It is simply a matter of verifying the three assertions below, which say that $\beta$ is obtained by following the procedure to construct form-reversing involutions on product spaces given in proposition \ref{anti-inv}.
\begin{lemma}\label{pdt_of_conj_classes}
Consider an integer $l\geq 1$ and let $\calC_1,\, ...\, ,\calC_l$ be $l$
conjugacy classes in $U$. Let $\beta^{(1)}$ be the map defined on
$\calC_1$ by
\begin{eqnarray*}
\beta^{(1)}: \calC_1 & \longto & \calC_1 \\
 c_1 & \longmapsto & \taum(c_1)
\end{eqnarray*}
and let $\beta^{(l-1)}$ be the map defined on the product
$\calC_2 \times\cdots\times \calC_l$ of $(l-1)$ conjugacy classes by:
\begin{eqnarray*}
\beta^{(l-1)}:  \calC_2 \times \cdots \times \calC_l & \longto &  \calC_2 \times \cdots \times \calC_l \\
(c_2,\, ...\, , c_l) & \longmapsto & (\taum(c_3...c_l)\taum(c_2)\tau(c_3..c_l),\, ...\, ,\taum(c_l))
\end{eqnarray*}
Let $\mu^{(l-1)}$ be the map from $\calC_2 \times\cdots\times \calC_l$ to $U$
defined by:
\begin{eqnarray*}
\mu^{(l-1)}: \calC_2 \times \cdots \times \calC_l & \longto &  U\\
(c_2,\, ...\, , c_l) & \longmapsto & c_2...c_l
\end{eqnarray*}
Finally, let $\beta^{(l)}$ be the map defined on
$\pconj=\calC_1\times (\calC_2 \times\cdots\times \calC_l)$ by:
\begin{eqnarray*}
\beta^{(l)}: \pconj & \longto & \pconj \\
(c_1,\, ...\, , c_l) & \longmapsto & (\taum(c_2...c_l)\taum(c_1)\tau(c_2...c_l),\, ...\ ,
\taum(c_l))
\end{eqnarray*}
Then we have:
$$\beta^{(l)}(c_1,\, ...\, ,c_l) = \Big(
\big(\mu^{(l-1)}\circ\beta^{(l-1)}(c_2,\, ...\,
,c_l)\big).\beta^{(1)}(c_1),\, \beta^{(l-1)}(c_2,\, ...\, ,c_l)\Big)$$
which we will write: $$\beta^{(l)}=
\Big(\big(\mu^{(l-1)}\circ\beta^{(l-1)}\big)
.\beta^{(1)},\, \beta^{(l-1)} \Big).$$
\end{lemma}
\begin{lemma}\label{pdt_of_doubles}
Consider an integer $g\geq 1$. Let $\beta_{(1)}$ be the map defined on
$U\times U$ by
\begin{eqnarray*}
\beta_{(1)}: U\times U & \longto & U\times U \\
 (a,b) & \longmapsto & (\tau(b),\tau(a))
\end{eqnarray*}
and let $\beta_{(g-1)}$ be the map defined on the product
$(U\times U)^{g-1}$ of $(g-1)$ copies of $(U\times U)$ by:
\begin{eqnarray*}
\beta_{(g-1)}:  (U\times U)^{g-1} & \longto &  (U\times U)^{g-1}\\
(a_2,b_2\, ...\, , a_g,b_g) & \longmapsto & (\taum([a_3,b_3]...[a_g,b_g])\tau(b_2)\tau([a_3,b_3]...[a_g,b_g]),\\
 & & \, ...\, ,\tau(a_g))
\end{eqnarray*}
Let $\mu_{(g-1)}$ be the map from $(U\times U)^{g-1}$ to $U$
defined by:
\begin{eqnarray*}
\mu_{(g-1)}: (U\times U)^{g-1} & \longto & U\\
(a_2, b_2\, ...\, , a_g,b_g) & \longmapsto & [a_2,b_2]...[a_g,b_g]
\end{eqnarray*}
Finally, let $\beta_{(g)}$ be the map defined on
$(U\times U)^g$ by:
\begin{eqnarray*}
\beta_{(g)}: (U\times U)^g & \longto & (U\times U)^g \\
(a_1, b_1\, ...\, , a_g,b_g) & \longmapsto & (\taum([a_2,b_2]...[a_g,b_g])\tau(b_1)\tau([a_2,b_2]...[a_g,b_g]),\\
& & \, ...\ ,
\tau(a_g))
\end{eqnarray*}
Then we have:
$$\beta_{(g)}(a_1, ... ,b_g) = \Big(
\big(\mu_{(g-1)}\circ\beta_{g-1)}(a_2, ...
,b_g)\big).\beta_{(1)}(a_1,b_1), \beta_{(g-1)}(a_2, ... ,b_g)\Big)$$
which we will write: $$\beta_{(g)}=
\Big(\big(\mu_{(g-1)}\circ\beta_{(g-1)}\big)
.\beta_{(1)},\, \beta_{(g-1)} \Big)$$
\end{lemma}
\begin{lemma}\label{ext_moduli_space}
Let $\beta^{(l)}$ and $\beta_{(g)}$ be defined as in lemmas \ref{pdt_of_conj_classes} and \ref{pdt_of_doubles} on $\pconj$ and $(U\times U)^g$ respectively. Let $\mu^{(l)}$ be the map from $\pconj$ to $U$ defined by $$\mu^{(l)}(c_1,\, ...\, ,c_l)=c_1...c_l.$$ Then the map $\beta$ on $\Mtot$ defined in definition \ref{def_beta} satisfies: $$\beta=\Big( (\mu^{(l)}\circ\beta^{(l)}).\beta_{(g)}, \beta^{(l)}\Big).$$
\end{lemma}
The above three lemmas are verified immediately. As a consequence, we have:
\begin{proposition}\label{beta_anti-inv}
The map $\beta$ on $\Mtot$ defined in \ref{def_beta} is an involution which reverses the $2$-form $\w$ defining the quasi-Hamiltonian structure on $\Mtot$ and which is compatible with the action of $U$ on this space and the momentum map $$\mu\abc=\relabcshort.$$ Additionally, $\beta$ has a non-empty fixed-point set $Fix(\beta)\not=\emptyset$, and $\mu(Fix(\beta))$ intersects the connected component $Q_0$ of $Fix(\taum)$ containing $1_U$: $$\mu(Fix(\beta))\cap Q_0\not=\emptyset.$$
\end{proposition}
\begin{proof}
Observe that $\beta^{(1)}$ is an involution that reverses the $2$-form defining the quasi-Hamiltonian structure on a single conjugacy class of $\calC$ of $U$ (see lemma \ref{anti-inv_on_conj_class}) and is compatible with the conjugation action of $U$ on $\calC$ and the momentum map $\mu^{(1)}:\calC\hookrightarrow U$. Likewise, $\beta_{(1)}$ is an involution that reverses the $2$-form defining the quasi-Hamiltonian structure on the double $U\times U$ (see lemma \ref{anti-inv_on_double}) and is compatible with the conjugation action of $U$ on $U\times U$ and the momentum map $\mu_{(1)}:(a,b)\in U\times U\mapsto aba^{-1}b^{-1}\in U$. Lemmas \ref{pdt_of_conj_classes} and \ref{pdt_of_doubles} show that we can apply proposition \ref{anti-inv} and therefore prove that $\beta^{(l)}$ and $\beta_{(g)}$ are form-reversing involutions on the quasi-Hamiltonian spaces $\pconj$ and $(U\times U)^g$ respectively, compatible with the respective actions of $U$ and their respective momentum maps. Lemma \ref{ext_moduli_space} then shows that we can apply proposition \ref{anti-inv} to $\beta$ and hereby prove that it is a form-reversing involution on the quasi-Hamiltonian space $\Mtot$ which is compatible with the action of $U$ and the momentum map $\mu\abc=\relabcshort$.\\
As for the fact that $\beta$ has fixed points, consider a maximal torus $T\subset U$ fixed pointwise by $\taum$. Such a torus is also globally $\tau$-stable. Take now an element $$(r_1,\tau(r_1),\, ...\, ,r_g,\tau(r_g),t_1,\, ...\, ,t_l)\in \big(\Mtot\big)\cap T^{2g+l}.$$ Since $T$ is Abelian and fixed pointwise by $\taum$, one obtains: $$\beta(r_1,\tau(r_1),\, ...\, ,r_g,\tau(r_g),t_1,\, ...\, ,t_l)=(r_1,\tau(r_1),\, ...\, ,r_g,\tau(r_g),t_1,\, ...\, ,t_l).$$ Finally, $$\mu(r_1,\tau(r_1),\, ...\, ,r_g,\tau(r_g),t_1,\, ...\, ,t_l)$$ is a product of elements of $T$ (namely, it is equal to $t_1...t_l$, all the commutators being equal to $1$ in the Abelian group $T$), hence $$\mu(r_1,\tau(r_1),\, ...\, ,r_g,\tau(r_g),t_1,\, ...\, ,t_l)\in T.$$ Since $T\subset Fix(\taum)$ is connected and contains $1_U$, we have $T\subset Q_0\subset Fix(\taum)$. Consequently $\mu(r_1,\tau(r_1),\, ...\, ,r_g,\tau(r_g),t_1,\, ...\, ,t_l)\in Q_0$, which ends the proof of the proposition.
\end{proof}
\begin{remark}\label{max_rk_assump_I}
As seen in the course of the proof of proposition \ref{beta_anti-inv}, the fact that $\beta$ has fixed points relies strongly on the fact that there exists a \emph{maximal} torus of $U$ (in particular, this torus contains an element of each conjugacy class) fixed \emph{pointwise} by $\taum$, which is the assumption made in section \ref{decomp}.
\end{remark}

\subsection{The characterization of decomposable representations of $\pi_{g,l}$}\label{charac_subsection}

We can now state our characterization result of decomposable representations of $\pi_{g,l}$:

\begin{theorem}[Characterization of decomposable representations]\label{charac_thm} \ 
Consider $$\abc\in \Mext$$ such that $\relabc=1$. Let $\beta$ be the involution introduced in definition \ref{def_beta}. Then the representation of $\pi_{g,l}=\pigl$ corresponding to $\abc$ is decomposable in the sense of definition \ref{decomp_rep} if and only if there exists $\pphi \in U$ such that $$\beta\abc=\pphi.\abc$$and$$\pphi\in Fix(\taum)\subset U.$$
\end{theorem}
\begin{proof}
We will denote $\abc$ by $(a,b,c)$.\\
\indent Let us first prove that if $(a,b,c)$ is a decomposable representation, then $\beta(a,b,c)$ $=\pphi.(a,b,c)$ for some $\pphi\in Fix(\taum)$. If $(a,b,c)$ is decomposable, there exists $(g+l)$ elements $v_1,\, ...\, ,v_g,w_1,\, ...\, ,w_l\in Fix(\taum)$ such that $[a_i,b_i]=v_i v_{i+1}^{-1}$, $c_j=w_jw_{j+1}^{-1}$ (with $v_{g+1}=w_1$ and $w_{l+1}=v_1$) and $\tau(a_i)=v_{i+1}^{-1}b_iv_{i+1}$. 
Therefore, we have, for all $j$ in $\{1,\, ...\, ,l\}$:
\begin{eqnarray*}
\taum(c_{j+1}...c_l)\taum(c_j)\tau(c_{j+1}...c_l) & = & \taum(w_{j+1}v_1^{-1}) \taum(w_jw_{j+1}^{-1}) \tau(w_{j+1} v_1^{-1}) \\
& = & v_1^{-1}w_{j+1} (w_{j+1}^{-1}w_j)w_{j+1}^{-1}v_1\\
& = & v_1^{-1} c_j v_1
\end{eqnarray*}
Also:
\begin{eqnarray*}
& & \taum([a_{i+1},b_{i+1}]...c_l) \tau(a_i) \tau([a_{i+1},b_{i+1}]...c_l) \\
& = & \taum(v_{i+1}v_1^{-1}) (v_{i+1}^{-1}b_i v_{i+1}) \tau(v_{i+1}v_1^{-1})\\
& = & v_1^{-1} v_{i+1} v_{i+1}^{-1} b_i v_{i+1} v_{i+1}^{-1} v_1\\
& = & v_1^{-1} b_i v_1
\end{eqnarray*}
And:
\begin{eqnarray*}
& & \taum([a_{i+1},b_{i+1}]...c_l) \tau(b_i) \tau([a_{i+1},b_{i+1}]...c_l) \\
& = & \taum(v_{i+1}v_1^{-1}) \tau(v_{i+1}\tau(a_i) v_{i+1}^{-1}) \tau(v_{i+1}v_1^{-1})\\
& = & v_1^{-1} v_{i+1} v_{i+1}^{-1} a_i v_{i+1} v_{i+1}^{-1} v_1\\
& = & v_1^{-1} a_i v_1
\end{eqnarray*}
Therefore, we have indeed $\beta(a,b,c)=v_1^{-1}.(a,b,c)$ with $v_1^{-1}\in Fix(\taum)$.\\
\indent Conversely, let us now assume that $\beta(a,b,c)=\pphi.(a,b,c)$ for some $\pphi\in Fix(\taum)$ and let us show that $(a,b,c)$ is decomposable. Since $\beta(a,b,c)=\pphi.(a,b,c)$, we have in particular: $\taum(c_l)=\pphi c_l\pphi^{-1}$, hence $\taum(c_l\pphi^{-1})=c_l\pphi^{-1}$, since $\taum(\pphi)=\pphi$. Therefore $c_l$ decomposes into a product of two symmetric elements: $c_l=(c_l\pphi^{-1})\pphi$. By induction, we obtain, using the equations $\beta(a,b,c)=\pphi.(a,b,c)$: 
$$\taum(c_{i+1}...c_l)\taum(c_i)\taum(c_{i+1}...c_l)=\pphi_i\pphi^{-1}=\tau(\pphi^{-1})c_i\taum(\pphi^{-1})$$ hence
\begin{eqnarray*}
\taum(c_ic_{i+1}...c_l\pphi^{-1}) & = & c_i\taum(\pphi^{-1})\taum(c_{i+1}...c_l)\\
& = & c_i \taum(c_{i+1}...c_l\pphi^{-1})\\
& = & c_i c_{i+1}...c_l\pphi^{-1}
\end{eqnarray*}
From the same set of equations $\beta(a,b,c)=\pphi.(a,b,c)$, we also have:
\begin{eqnarray}
\taum([a_{i+1},b_{i+1}]...c_l) \tau(a_i) \tau([a_{i+1},b_{i+1}]...c_l) & = & \pphi b_i \pphi^{-1} \label{eq1}\\
\taum([a_{i+1},b_{i+1}]...c_l) \tau(b_i) \tau([a_{i+1},b_{i+1}]...c_l) & = & \pphi a_i \pphi^{-1} \label{eq2}
\end{eqnarray}
Taking the commutator of (\ref{eq2}) and (\ref{eq1}), we obtain: $$\pphi[a_i,b_i]\pphi^{-1} = \taum([a_{i+1},b_{i+1}]...c_l)\underbrace{[\tau(b_i),\tau(a_i)]}_{=\taum([a_i,b_i])}\tau([a_{i+1},b_{i+1}]...c_l)$$
hence: $$[a_i,b_i]\pphi^{-1}\taum([a_{i+1},b_{i+1}]...c_l)=\pphi^{-1}\taum([a_i,b_i]...c_l).$$
Using the induction hypothesis $$\taum([a_{i+1},b_{i+1}])...c_l\pphi^{-1})=[a_{i+1},b_{i+1}]...c_l\pphi^{-1}$$ and the fact that $\taum(\pphi)=\pphi$, we obtain: $$\taum([a_i,b_i]...c_l\pphi^{-1})=[a_i,b_i]...c_l\pphi^{-1}.$$ We can then set:
\begin{eqnarray*}
v_i & := & [a_i,b_i]...c_l\pphi^{-1} \in Fix(\taum)\\
w_j & := & c_j...c_l\pphi^{-1} \in Fix(\taum)
\end{eqnarray*}
Since $(a,b,c)$ is a representation, we have $\relabc=1$, hence $v_1=\pphi^{-1}$. Therefore, we have: $$[a_i,b_i]=v_iv_{i+1}^{-1}\ \mathrm{and}\ c_j=w_jw_{j+1}^{-1}$$ (with $v_{g+1}:=w_1$ and $w_{l+1}:=v_1$). Finally, we can then rewrite equation (\ref{eq1}) under the form: $$\taum(v_{i+1}v_1^{-1})\tau(a_i)\tau(v_{i+1}v_1^{-1}) = \pphi b_i \pphi^{-1} = v_1^{-1} b_i v_1$$ hence:
\begin{eqnarray*}
\tau(a_i) & = & \tau(v_{i+1}v_1^{-1})v_1^{-1} b_i v_1 \taum(v_{i+1}v_1^{-1})\\
& = & v_{i+1}^{-1} v_1 v_1^{-1} b_i v_1 v_1^{-1} v_{i+1}\\
& = & v_{i+1}^{-1}  b_i v_{i+1}
\end{eqnarray*}
This concludes the proof that the representation $(a,b,c)$ is decomposable in the sense of definition \ref{decomp_rep}.
\end{proof}
\begin{remark}\label{nice_case}
If $Fix(\taum)$ is connected, which is the case for instance for $\tau(u)=\overline{u}$ on $U(n)$, a symmetric element $\pphi\in Fix(\taum)$ may be written $\pphi=\taum(u)u$ for some $u\in U$. In this case, the condition $\beta(a,b,c)=\pphi.(a,b,c)$ with $\pphi=\taum(u)u\in Fix(\taum)$ implies $\beta(u.(a,b,c))=u.(a,b,c)$. This is the point of view adopted in \cite{Sch_CJM}.
\end{remark}
We now know that the decomposable representations of $\pi_{g,l}$ are exactly the representations $(a,b,c)$ satisfying $\beta(a,b,c)=\pphi.(a,b,c)$ for some $\pphi\in Fix(\taum)$. In the following section, we will use this characterization to prove that decomposable representations always exist. In fact, we will prove a stronger result: there always exist representations of $\pi_{g,l}$ which are \emph{fixed} by $\beta$, that is, $Fix(\beta)\cap\mu^{-1}(\{1\})\not=\emptyset$.

\section{Existence of decomposable representations}\label{existence}

In this section, we will show that there always exist decomposable representations of the surface group $\pi_{g,l}$ into an arbitrary compact connected Lie group $U$. This proof of existence relies on their characterization in terms of an involution $\beta$ obtained in theorem \ref{charac_thm} and on a real convexity theorem for group-valued momentum maps (\cite{Sch_IM}).\\ Recall from proposition \ref{beta_anti-inv} that $\beta$ is a form-reversing involution on the quasi-Hamiltonian space $$M=\Mtot$$ compatible with the diagonal conjugation action of $U$ and with the momentum map $$\mu\abc=\relabc$$ of this action. Consequently, the involution $\beta$ induces an anti-symplectic involution $\bhat$ on the representation space $$\Mgl=\qhamquot$$ defined by $\bhat([a,b,c]):=[\beta(a,b,c)]$. As announced at the end of section \ref{charac}, we will now prove the following result: $$Fix(\beta)\cap\mu^{-1}(\{1\})\not=\emptyset.$$ In view of theorem \ref{charac_thm}, this shows that there exist decomposable representations (which are representations $(a,b,c)$ satisfying $\beta(a,b,c)=\pphi.(a,b,c)$ with $\pphi\in Fix(\taum)$) since in particular $1\in Fix(\taum)$). Observe that our choice of conjugacy classes $\calC_1,\, ...\, ,\calC_l\subset U$ is assumed to be such that $\mu^{-1}(\{1\})\not=\emptyset$ (that is: the representation space is not empty) and that we know from proposition \ref{beta_anti-inv} that $Fix(\beta)\not=\emptyset$. In order to prove that $Fix(\beta)\cap\mu^{-1}(\{1\})\not=\emptyset$, we will distinguish three cases. We will begin with the case where the compact connected Lie group $U$ is in addition simply connected. In this case, the fact that $Fix(\beta)\cap\mu^{-1}(\{1\})\not=\emptyset$ is a corollary of the real convexity theorem for group-valued momentum maps obtained in \cite{Sch_IM}. We will then deal with the case of compact connected semi-simple Lie groups and eventually move on to arbitrary compact connected Lie groups. In these last two cases, we will reduce the situation to the case of simply connected groups (using covering spaces).

\subsection{The case where $U$ is simply connected}\label{simply_connected_case}

When $U$ is a compact connected simply connected Lie group, the existence of decomposable representations of the surface group $$\pi_{g,l}=\pigl$$ is a consequence of the following convexity result for images of fixed-point sets of form-reversing involutions defined on quasi-Hamiltonian spaces (see \cite{OSS} for a real convexity theorem in the usual Hamiltonian case, which the following statement is a quasi-Hamiltonian analogue of).
\begin{theorem}[A real convexity theorem for group-valued momentum maps, \cite{Sch_IM}]\label{convexity_thm}
Let $(U,\tau)$ be a compact connected simply connected Lie group endowed with an involutive automorphism $\tau$ such that the involution $\taum:u\in U\mapsto \tau(u^{-1})$ leaves a maximal torus $T$ of $U$ pointwise fixed, and let $\Walc\subset \mathfrak{t}:=Lie(T)$ be a Weyl alcove.
Let $(M,\w,\mu:M\to U)$ be a connected quasi-Hamiltonian $U$-space with proper momentum map
$\mu:M\to U$ and let $\beta:M\to M$ be an involution on $M$ satisfying:
\begin{enumerate}
\item[(i)] $\beta^*\w=-\w$.
\item[(ii)] $\beta(u.x)=\tau(u).\beta(x)$ for all $x\in M$ and all $u\in U$.
\item[(iii)] $\mu\circ\beta=\taum\circ\mu$.
\item[(iv)] $M^{\beta}:=Fix(\beta)\not=\emptyset$.
\item[(v)] $\mu(M^{\beta})$ has a non-empty intersection with the fixed-point set $Q_0$ of $1$ in $Fix(\taum)\subset U$.
\end{enumerate}
\noindent Then: $$\mu(M^{\beta})\cap\exp(\clWalc) =\mu(M)\cap\exp(\clWalc).$$
In particular, $\mu(M^{\beta})\cap\exp(\clWalc)$ is a convex subpolytope of $\exp(\clWalc)\simeq\clWalc\subset\mathfrak{t}$, equal to the full momentum polytope $\mu(M)\cap\exp(\clWalc)$.
\end{theorem}
\begin{remark}\label{max_rk_assump_II}
This theorem, which is the key to proving the existence of decomposable representations, uses the assumption that the involution $\tau$ is of maximal rank in a crucial way (see \cite{Sch_IM} for a detailed proof).
\end{remark}
As announced in remark \ref{fix_bhat}, we point out the following corollary (the semi-simple case will be treated in theorem \ref{inv_semi-simple}):
\begin{corollary}[$Fix(\beta)\cap\mu^{-1}(\{1\})\not=\emptyset$]\label{existence_cor}
If $\beta$ satisfies the assumptions of theorem \ref{convexity_thm} and $\bhat$ designates the induced involution $\bhat([x]):=[\beta(x)]$ on the quasi-Hamiltonian quotient $M//U=\qhamquot$, we have: $Fix(\beta)\cap\mu^{-1}(\{1\})\not=\emptyset$ and therefore $Fix(\bhat)\not=\emptyset$.
\end{corollary}
\begin{proof}
Since $\mu^{-1}(\{1\})\not=\emptyset$ and since we always have $1\in\exp(\clWalc)$, we obtain, using theorem \ref{convexity_thm}: $$1\in\mu(M)\cap\exp(\clWalc)=\mu(M^{\beta})\cap\exp(\clWalc)$$ that is: $$Fix(\beta)\cap\mu^{-1}(\{1\})\not=\emptyset.$$ If $x \in Fix(\beta)\cap\mu^{-1}(\{1\})\not=\emptyset$, then by definition $\bhat([x])=[\beta(x)]=[x]$.
\end{proof}
We now specify this to the case where $M$ is the quasi-Hamiltonian space $M=\Mtot$ and $\beta$ is the involution defined in \ref{def_beta}.
\begin{theorem}[Existence of decomposable representations, (I)]\label{existence_thm_I}
If $U$ is a compact connected simply connected Lie group, there always exist decomposable representations of $\pi_{g,l}=\pigl$ into $U$.
\end{theorem}
\begin{proof}
We know from theorem \ref{charac_thm} that there exist decomposable representations if and only if there exist representations $(a,b,c)$ of $\pi_{g,l}$ satisfying $\beta(a,b,c)=\pphi.(a,b,c)$ where $\beta$ is the involution defined in definition \ref{def_beta}. Proposition \ref{beta_anti-inv} shows that $\beta$ satisfies the assumptions of theorem \ref{convexity_thm} (in particular assumption (v) is satisfied), hence corollary \ref{existence_cor} applies and shows that we always have $Fix(\beta)\cap\fiber\not=\emptyset$, which proves the theorem since in that case $\beta(a,b,c)=(a,b,c)=1.(a,b,c)$ and $1\in Fix(\taum)$.
\end{proof}

\subsection{The case where $U$ is semi-simple}\label{semi-simple_case}

When $U$ is a compact connected semi-simple Lie group, its universal cover is still compact. Let us denote it by $\tU$, and by $\pi$ the covering map $\pi:\tU\to U$. We will then borrow from \cite{AMW2} to construct a quasi-Hamiltonian $\tU$-space associated to our $M=\Mtot$. I would like to thank Eckhard Meinrenken for bringing to my attention the reference \cite{AMW2}.\\
Following \cite{AMW2}, we set: $$\tM:=M\times_U\tU =\{(x,\tu)\ |\ \mu(x)=\pi(\tu)\}$$ and
\begin{eqnarray*}
p:\tM \longto M & & \tmu:\tM\longto \tU \\
(x,\tu) \longmapsto x & & (x,\tu) \longmapsto \tu
\end{eqnarray*}
so that we have the following commutative diagram:
$$\begin{CD}
\tM=M\times_U \tU @>\tmu>> \tU\\
@VpVV  @VV\pi V \\
M @>\mu>> U
\end{CD}$$
Finally, let us set $\tw:=p^*\w$ and observe that $\tU$ acts on $\tM$ via $$\tu_0.(x,\tu):=(\pi(\tu_0).x,\tu_0\tu\tu_0^{-1})$$ and that $\tmu$ is equivariant for this action. We then quote from \cite{AMW2}:
\begin{proposition}[\cite{AMW2}]\label{covering}
$(\tM=M\times_U \tU,\tw,\tmu:\tM\to\tU)$ is a quasi-Hamiltonian $\tU$-space.
\end{proposition}
The proof shows that this works because $\pi:\tU\to\tU$ is a \emph{covering} homomorphism. In particular, $\tU$ and $U$ have isomorphic Lie algebras. Further, we have:
\begin{proposition}[\cite{AMW2}]\label{covering_quotient}
The quasi-Hamiltonian quotients associated to $M$ and $\tM$ are isomorphic: the map $p:\tM\to M$ sends $\tmu^{-1}(\{1_{\tU}\})$ to $\mu^{-1}(\{1_U\})$ and induces an isomorphism $$\tmu^{-1}(\{1_{\tU}\})/\tU\simeq\mu^{-1}(\{1_U\})/U.$$ In particular, if $\mu^{-1}(\{1_U\})\not=\emptyset$ then $\tmu^{-1}(\{1_{\tU}\})\not=\emptyset$.
\end{proposition}
Let us now observe that $\beta$ induces a form-reversing involution $\tbeta$ on $\tM$.
\begin{lemma}\label{group_inv}
Since the compact connected groups $U$ and $\tU$ have isomorphic Lie algebras, the involutive automorphism $\tau$ of $U$ induces an involutive automorphism of $\tU$, that we denote by $\ttau$. In particular, we have $\pi\circ\ttau=\tau\circ\pi$, where $\pi$ is the covering map $\pi:\tU\to U$. As usual, we will denote by $\ttaum$ the involution $\ttaum(\tu):=\ttau(\tu^{-1})$. If $\tau$ is of maximal rank, so is $\ttau$. If we denote by $Q_0$ the connected component of $1_U$ in $Fix(\taum)\subset U$ and by $\tQ_0$ the connected component of $1_{\tU}$ in $Fix(\tau)\subset \tU$, the covering map $\pi:\tU\to U$ restricts to a covering map $\pi|_{\tQ_0}:\tQ_0\to Q_0$.
\end{lemma}
\begin{proposition}\label{covering_inv}
Let $\beta$ be a form-reversing involution on the quasi-Hamiltonian space $(M,\w,\mu:M\to U)$, compatible with the action of $(U,\tau)$ and the momentum map $\mu$. Then the map
\begin{eqnarray*}
\tbeta: \tM & \longto & \tM \\
(x,\tu) & \longmapsto & (\beta(x),\ttaum(\tu))
\end{eqnarray*}
is a form-reversing involution on the quasi-Hamiltonian space $(\tM=M\times_U \tU,\tw,\tmu:\tM\to\tU)$, compatible with the action of $(\tU,\ttau)$ and the momentum map $\tmu$.
\end{proposition}
\begin{remark}\label{on_coverings}
Observe that propositions \ref{covering}, \ref{covering_quotient} and \ref{covering_inv}, as well as lemma \ref{group_inv},  are true for any covering map $\pi:\tU\to U$, regardless of whether $\tU$ is the universal cover of $U$. This will be useful in subsection \ref{arbitrary_case}.
\end{remark}
We then have:
\begin{theorem}\label{inv_semi-simple}
Let $(U,\tau)$ be a compact connected semi-simple Lie group endowed with an involutive automorphism $\tau$ of maximal rank, and let $(M,\w,\mu:M\to U)$ be a connected quasi-Hamiltonian $U$-space such that $\mu^{-1}(\{1_U\})\not=\emptyset$.  Let $\beta$ be a form-reversing compatible involution $\beta$ on $M$, whose fixed-point set $Fix(\beta)$ is not empty and has an image under $\mu$ that intersects the connected component of $1_U$ in $Fix(\taum)\subset U$. Then: $$Fix(\beta)\cap\mu^{-1}(\{1_U\})\not=\emptyset.$$
\end{theorem}
\begin{proof}
We will show that there exists a connected component of $\tM=M\times_U\tU$ which contains points of $\tmu^{-1}(\{1_{\tU}\})$ \emph{and} fixed points of $\tbeta$, and apply the corollary of the convexity theorem (corollary \ref{existence_cor}) to this connected component, which is a quasi-Hamiltonian space. From this we will deduce the statement of the theorem.\\
Since $\mu^{-1}(\{1\})\not=\emptyset$ and $\mu(Fix(\beta))\cap Q_0\not=\emptyset$, there exist $x_0\in M$ such that $\mu(x_0)=1_U$ and $x_1\in M$ such that $\beta(x_1)=x_1$ and $\mu(x_1)\in Q_0$. Since $M$ is connected, there is a path $(x_t)_{t\in [0,1]}$ from $x_0$ to $x_1$. Set $u_t:=\mu(x_t)\in U$ for all $t\in [0,1]$. Since $\pi:\tU\to U$ is a covering map, we can lift the path $(u_t)_{t\in [0,1]}$ to a path $(\tu_t)_{t\in [0,1]}$ on $\tU$ such that $\pi(\tu_t)=u_t=\mu(x_t)$ and $\tu_0=1_{\tU}$. Then $(x_t,\tu_t)\in \tM=M\times_U \tU$ and it is a path going from $(x_0,\tu_0)=(x_0,1_{\tU})$ to $(x_1,\tu_1)$, which are therefore contained in a same connected component $\tM_0$ of $\tM$. Then, we have $\tmu(x_0,1_{\tU})=1_{\tU}$ and, since $\pi(\tu_1)=u_1=\mu(x_1)\in Q_0\subset Fix(\taum)$, we have $\tu_1\in\tQ_0\subset Fix(\ttaum)$, hence $$\tbeta(x_1,\tu_1)=(\beta(x_1),\ttaum(\tu_1))=(x_1,\tu_1)$$ and $\tmu(x_1,\tu_1)=\tu_1\in\tQ_0$. Therefore, the connected component $\tM_0$ of $\tM$, which is a quasi-Hamiltonian $\tU$-space, contains points of $\tmu^{-1}(\{1_{\tU}\})$ and points of $Fix(\tbeta)$ whose image is contained in $\tQ_0$. Since $\tU$ is simply connected, we can apply corollary \ref{existence_cor} and conclude that $Fix(\tbeta)\cap\tmu^{-1}(\{1_{\tU}\})\not=\emptyset$. Take now $(x,\tu)\in Fix(\tbeta)\cap\tmu^{-1}(\{1_{\tU}\})$. In particular, $\tu=1_{\tU}$. Since $\tbeta(x,\tu)=(x,\tu)$, we have $\beta(x)=x$ and $\mu(x)=\mu\circ p(x,\tu)=\pi\circ\tmu(x,\tu)=\pi(\tu)=\pi(1_{\tU})=1_U$. That is: $x\in Fix(\beta)\cap\mu^{-1}(\{1_U\})$, which is therefore non-empty.
\end{proof}
In view of the characterization of decomposable representations obtained in theorem \ref{charac_thm}, we then have:
\begin{theorem}[Existence of decomposable representations, (II)]\label{existence_thm_II}
If $U$ is a compact connected semi-simple Lie group, there always exist decomposable representations of $\pi_{g,l}=\pigl$ into $U$.
\end{theorem}
\begin{proof}
This is proved similarly to theorem \ref{existence_thm_I}: the proof follows from applying theorem \ref{inv_semi-simple} to the quasi-Hamiltonian space $M=\Mtot$ and from the characterization of decomposable representations obtained in theorem \ref{charac_thm}.
\end{proof}

\subsection{The case of an arbitrary compact connected Lie group}\label{arbitrary_case}

To handle the case of an arbitrary compact connected Lie group $U$, we will use the same technique of coverings as in subsection \ref{semi-simple_case}, but instead of considering the universal cover (which is not compact if $U$ is not semi-simple), we will use the following result on the structure of compact Lie groups: any compact connected Lie group $U$ possesses a finite cover $\tU$ of the form $\tU=S\times G$, where $S\subset \Z(U)$ is contained in the center of $U$ (in fact, it is the neutral component of $\Z(U)$, in particular, it is an Abelian Lie group) and $G\subset U$ is a connected Lie group which is both compact \emph{and} simply connected (see for instance \cite{BtD}) and whose Lie algebra is $\mathfrak{g}=[\mathfrak{u},\mathfrak{u}]$ with covering map $\pi:(s,g)\in S\times G\mapsto sg\in U$. Observe here that $\tU$ is \emph{not} the universal cover of $U$ unless $U$ is semi-simple. The same type of approach to arbitrary compact groups is used in \cite{AMW2} and \cite{Slee}. I would like to thank Pierre Sleewaegen for discussions on this and Eckhard Meinrenken for pointing out to me the reference \cite{AMW2}.\\
If now $(M,\w,\mu:M\to U)$ is a quasi-Hamiltonian space, then by proposition \ref{covering}, we can construct (see remark \ref{on_coverings}) a quasi-Hamiltonian space $$(\tM,\tw,\tmu:\tM\to \tU)$$ this time with $\tU=S\times G$. Let us now briefly come back to the example we will eventually apply all this to: $$M=\Mtot$$ endowed with the diagonal conjugation action of $U$. Then $$\tM=M\times_U(S\times G)$$ is endowed with the following action of $S\times G$: $$(s_0,g_0).(x,s,g)=\big(\pi(s_0,g_0).x,s_0ss_0^{-1},g_0gg_0^{-1}\big)$$ where $\pi$ is the covering map $\pi:\tU=S\times G\to U$. Since $S$ is Abelian, this rewrites $(s_0ss_0^{-1},g_0gg_0^{-1})=(s,g_0gg_0^{-1})$. Further, when $$M=\Mtot$$ we observe that $S\subset\Z(U)$ acts trivially on $M$. Therefore we have here: $$\pi(s_0,g_0).x=s_0.(g_0.x)=g_0.x.$$ That is: $S$ \emph{acts trivially on} $\tM$. \emph{We will now work under the assumption that} $S$ \emph{acts trivially on} $\tM$. If we denote by $p_G$ the projection $p_G:S\times G\to G$ and by $\mu_G$ the map $\mu_G:=p_G\circ\tmu:\tM\to G$, we then have:
\begin{proposition}\label{projection}
If $S\hookrightarrow \tU=S\times G$ acts trivially on the quasi-Hamiltonian $\tU$-space $(\tM,\tw,\tmu:\tM\to\tU=S\times G)$, then $(\tM,\tw,\mu_G:\tM\to G)$ is a quasi-Hamiltonian $G$-space for the action of $G\hookrightarrow \tU=S\times G$.
\end{proposition}
\begin{proof}
This will be a consequence of the fact that $S$ is both Abelian and contained in the center of $U$ (hence also centralizes $G\subset [U,U]$).\\
Since $\tU=S\times G$ is a direct product, one can express the Maurer-Cartan $1$-forms and the Cartan $3$-form of $\tU$ in terms of those and $S$ and $G$. For instance, the Cartan $3$-form is $\chi_{\tU}=\chi_S\oplus\chi_G=0\oplus\chi_G$, since $S$ is Abelian. This proves that $$d\tw=-\tmu^*\chi_{\tU}=-\tmu(p_G^*\chi_G)=-(p_G\circ\tmu)^*\chi_G=-\mu_G^*\chi_G.$$ Further, we have, for all $y\in\tM$: $$\ker\tw_y=\{\tX^{\sharp}_y\ :\ \tX\in\mathfrak{u}=\mathfrak{s}\oplus\mathfrak{g}\ |\ Ad\, \tmu(y).\tX=-\tX\}.$$ But $\tX=(\xi,X)$ with $\xi\in\mathfrak{s}$ and $X\in\mathfrak{g}$. If we denote $\tmu(y)=(s,g)\in S\times G$, we have: $Ad\, \tmu(y).\tX=(Ad\, s.\xi,Ad\, g.X)=(\xi,Ad\, g.X)$, since $S$ is Abelian. Therefore, $Ad\, \tmu(y).\tX=-\tX=-(\xi,X)$ if and only if $\xi=0$ and $Ad\, g.X=-X$, that is: $Ad\, \mu_G(y).X=-X$. Now, since $S$ acts trivially on $\tM$, we have in fact $\tX^{\sharp}_y=X^{\sharp}_y$, where the action of $G$ is given by $G\hookrightarrow S\times G$. Indeed: $\tX^{\sharp}_y=\frac{d}{dt}|_{t=0}(\exp(t\tX).y)$ with $\tX=(\xi,X)$. Since $\xi$ centralizes X, we have $\exp(t\tX)=\exp(t\xi)\exp(tX)$, and since $S$ acts trivially on $\tM$, we have $\exp(t\xi)\exp(tX).y=\exp(tX).y$, therefore $\tX^{\sharp}_y=X^{\sharp}_y$. Consequently, the above proves that: $$\ker \tw_y=\{X^{\sharp}_y\ :\ X\in\mathfrak{g}\subset\mathfrak{u}=\mathfrak{s}\oplus\mathfrak{g}\ |\ Ad\, \mu_G(y).X=-X\}.$$ Finally, for all $X\in\mathfrak{g}$, setting $\tX=(0,X)\in\mathfrak{u}$, one has:
\begin{eqnarray*}
\iota_{X^{\sharp}}\tw & = & \iota_{\tX^{\sharp}}\tw\\
& = & \frac{1}{2} \tmu(\theta^L_{\tU}+\theta^R_{\tU}\ |\ (0,X))_{\widetilde{\mathfrak{u}}}\\
& = & (\theta_S\ |\ 0)_{\mathfrak{s}} + \frac{1}{2}(\theta^L_G+\theta^R_G\ |\ X)_{\mathfrak{g}}\\
& = & \frac{1}{2}(\theta^L_G+\theta^R_G\ |\ X)_{\mathfrak{g}}
\end{eqnarray*}
which proves that $(\tM,\tw,\mu_G:\tmu\to G)$ is a quasi-Hamiltonian $G$-space.
\end{proof}
Let us continue our analysis of the above situation, assuming that $S$ acts trivially on $\tM$. Recall from proposition \ref{covering_quotient} that $\tmu^{-1}(\{1_{\tU}\})\not=\emptyset$. The fact that $S$ acts trivially on $\tM$ then has the following consequence:
\begin{proposition}\label{constant_momentum}
Consider $(x_0,1_S,1_G)\in\tmu^{-1}(\{1_{\tU}=(1_S,1_G)\})$ and denote by $\tM_0$ the connected component of $\tM$ containing $(x_0,1_S,1_G)$. If $S\hookrightarrow S\times G=\tU$ acts trivially on $\tM$, then $\tmu(\tM_0)\subset \{1_S\}\times G$.
\end{proposition}
This is true because $\tmu$ is a momentum map: when the action is trivial, the momentum map is a constant map. Here, one factor of $S\times G$ acts trivially and consequently the component of the momentum map with values in this factor is a constant map.
\begin{proof}
Recall first that since $S$ centralizes $G$, the actions of $S$ and $G$ on $\tM_0$ commute. Therefore, the isotropy Lie algebra of an arbitary point $y\in \tM_0$ is $\widetilde{\mathfrak{u}}_y=\mathfrak{s}_y\oplus\mathfrak{g}_y$, where $\mathfrak{s}_y$ (resp. $\mathfrak{g}_y$) is the isotropy algebra for the action of $S$ (resp. $G$). Since $S$ acts trivially on $\tM_0$, we have in fact $\mathfrak{s}_y=\mathfrak{s}$. Then: $$\widetilde{\mathfrak{u}}_y^{\perp}=\{(\eta,Z)\in \mathfrak{s}\oplus\mathfrak{g}\ |\ \forall (\xi,X)\in\widetilde{\mathfrak{u}}_y=\mathfrak{s}_y\oplus\mathfrak{g}_y, (\eta\ |\ \xi)_{\mathfrak{s}} + (Z\ |\ X)_{\mathfrak{g}}=0\}.$$ Since $\mathfrak{s}_y=\mathfrak{s}$,we have: if $(\eta,Z)\in\widetilde{\mathfrak{u}}_y^{\perp}$, then for all $\xi\in\mathfrak{s}$, $(\eta\ |\ \xi)_{\mathfrak{s}}=0$, hence $\eta=0$. Therefore: $$\widetilde{\mathfrak{u}}_y^{\perp}=\{0\}\oplus\{ Z\in\mathfrak{g}\ |\ \forall X\in\mathfrak{g}_y, (Z\ |\ X)=0\}\simeq\mathfrak{g}_y^{\perp_{\mathfrak{g}}}\subset \mathfrak{g}.$$ But since $\tmu$ is the momentum map of the action of $\tU$ on $\tM_0$, we have (see \cite{AMM}) that $\im\, T_y\tmu\simeq \widetilde{\mathfrak{u}}_y^{\perp}$. Consequently: $\im\, T_y\tmu \simeq \widetilde{\mathfrak{u}}_y^{\perp}= \{0\}\oplus\mathfrak{g}_y^{\perp_{\mathfrak{g}}}\subset \{0\}\oplus\mathfrak{g}$. Therefore, since $\tM_0$ is connected, we have $\tmu(\tM_0)\subset \{c\}\times G\subset \tU$ for some $c\in S$. But since $(x_0,1_S,1_G)\in \tM_0$, we have $\tmu(x_0,1_S,1_G)\in \{c\}\times G$, hence $c=1_S$ and $\tmu(\tM_0)\subset \{1_S\}\times G$.
\end{proof}
We are now in the following situation: $(\tM_0,\tw,\mu_G:\tM_0\to G)$ is a quasi-Hamiltonian $G$-space, with $G$ compact connected and simply connected, and with $\mu_G^{-1}(\{1_G\})\not=\emptyset$. Further, as in proposition \ref{covering_inv}, the involution $\beta$ induces an involution $\tbeta$ on $\tM_0$. More specifically:
\begin{proposition}\label{fixed_points_of_tbeta}
The involution $\tbeta$ on $\tM$ induced by the involution $\beta$ on $M$ as in proposition \ref{covering_inv} satisfies $\tbeta(\tM_0)\subset\tM_0$ and $Fix(\tbeta|_{\tM_0})\not=\emptyset$. Further, the involution $\ttau$ on $\tU=S\times G$ restricts to an involution $\tau_G$ on $G$, which is of maximal rank if $\ttau$ is of maximal rank. The involution $\tbeta|_{\tM_0}$ is compatible with the action of $G$ on $\tM_0$ and the momentum map $\mu_G:\tM_0\to G$. Denoting by $Q_0^G$ the connected component of $Fix(\taum_G)$ containing $1_G$, then $\mu(Fix(\tbeta|_{\tM_0}))\cap Q_0^G\not=\emptyset$.
\end{proposition}
\begin{proof}
Observe that $\ttau(G)\subset G$ because $\mathfrak{g}=[\mathfrak{u},\mathfrak{u}]$. Now, as in the proof of theorem \ref{inv_semi-simple}, we start with a path $(x_t)_{t\in [0,1]}$ going from $x_0\in M$ (the \emph{same} $x_0$ used to define $\tM_0$ in proposition \ref{constant_momentum}) to some fixed point $x_1$ of $\beta$ satisfying $\mu(x_1)\in Q_0\subset U$. We can lift the path $u_t=\mu(x_t)\in U$ to a path $\tu_t\in \tU$ starting at $(1_S,1_G)\in\tU$. The path $(x_t,\tu_t)$ is then a path in $\tM$ starting at $(x_0,1_S,1_G)\in\tM_0$ and is therefore a path in $\tM_0$. Since $\tmu(\tM_0)\subset \{1_S\}\times G$, the path $(\tu_t)$ is of the form $\tu_t=(1_S,g_t)$.  Then, $(x_1,1_S,g_1)$ is a fixed point of $\tbeta$ which is \emph{contained in} $\tM_0$. Since $\tM_0$ is connected and contains a fixed point of $\tbeta$, we have $\tbeta(\tM_0)\subset\tM_0$, and, as in the proof of theorem \ref{inv_semi-simple}, we have $\mu_G(x_1,1_S,g_1)\in Q_0^G$.
\end{proof}
Since $G$ is simply connected, we can now apply the real convexity theorem \ref{convexity_thm} to $(\tM_0,\tw,\mu_G:\tM_0\to G)$ and we obtain:
\begin{theorem}\label{inv_arbitrary}
Let $(U,\tau)$ be an arbitrary compact connected Lie group endowed with an involutive automorphism of maximal rank $\tau$ and let $(M,\w,\mu:M\to U)$ be a connected quasi-Hamiltonian $U$-space such that $\mu^{-1}(\{1_U\})\not=\emptyset$. Let $\beta$ be a compatible form-reversing involution on $M$ whose fixed-point set $Fix(\beta)$ is non-empty and has an image under $\mu$ that intersects the connected component of $1_U$ in $Fix(\taum)\subset U$. Assume that the center $\Z(U)$ of $U$ acts trivially on $M$. Then: $$Fix(\beta)\cap\mu^{-1}(\{1_U\})\not=\emptyset.$$
\end{theorem}
\begin{proof}
Denote by $\tU=S\times G$ the finite cover of $U$ where $S\subset \Z(U)$ is a torus and $G\subset [U,U]$ is compact connected and simply connected. Since $S\subset \Z(U)$, $S$ acts trivially on $M$ and since it is in addition an Abelian group, we have seen that it also acts trivially on the lifted quasi-Hamiltonian space $\tM=M\times_U\tU$ introduced in proposition \ref{covering}. Propositions \ref{projection}, \ref{constant_momentum} and \ref{fixed_points_of_tbeta} then show that we have a connected quasi-Hamiltonian $G$-space $(\tM_0,\tw,\mu_G:\tM_0\to G)$ containing points of $\mu_G^{-1}(\{1_G\})$ and fixed points of $\tbeta$, whose image is in addition contained in the connected component of $1_G$ in $Fix(\taum_G)\subset G$ (where $\tau_G$ is the involution of maximal rank of $G$ induced by the involution $\tau$ of $U$). Moreover, we have $\tmu(\tM_0)\subset \{1_S\}\times G$. Since $G$ is compact connected and simply connected, we can apply corollary \ref{existence_cor} and conclude that $Fix(\tbeta|_{\tM_0})\cap\mu_G^{-1}(\{1_G\})\not=\emptyset$. Take now $(x,c,1_G)\in Fix(\beta|_{\tM_0})\cap\mu_G^{-1}(\{1_G\})$. Since $\tmu(M_0)\subset \{1_S\}\times G$, we have in particular, $c=1_S$. Consequently, $\mu(x)=\pi(1_S,1_G)=1_U$. And since $(x,1_S,1_G)\in Fix(\tbeta)$, we have in particular $\beta(x)=x$. Hence $x\in Fix(\beta)\cap\mu^{-1}(\{1_U\})$, which is therefore non-empty.
\end{proof}
As a consequence of all this, we have:
\begin{theorem}[Existence of decomposable representations, (III)]\label{existence_thm_III}
If $U$ is an arbitrary compact connected Lie group, there always exist decomposable representations of $\pi_{g,l}=\pigl$ into $U$.
\end{theorem}
\begin{proof}
The proof follows from applying theorem \ref{inv_arbitrary} to the quasi-Hamiltonian space $M=\Mtot$, on which the center $\Z(U)$ of $U$ acts trivially (since the action of $U$ on $M$ is the diagonal conjugation action) and from the characterization of decomposable representations obtained in theorem \ref{charac_thm}.
\end{proof}
We now have all the ingredients that we need to obtain Lagrangian submanifolds of the representation spaces $\Mgl=\Mod$. The fact that there always exist decomposable representations will be used under the form $Fix(\bhat)\not=\emptyset$, where $\bhat$ is the involution induced by $\beta:M\to M$ on $\Mgl=M//U$. Let us mention here that the results of this section, specifically theorem \ref{existence_thm_III}, give an alternative proof of the existence of decomposable representations of $\piS$ into $U(n)$, a result first obtained by Falbel and Wentworth in \cite{FW}, by methods of a completely different nature.  

\section{The Lagrangian nature of decomposable representations}\label{lag_nature}

In this section we will achieve the goal of this paper, which is to give an example of Lagrangian submanifold of the representation space $$\Mgl=\Mod=\Mtot//U$$ associated to the surface group $$\pi_{g,l}=\pigl$$ where $\Sigma_g$ is a compact Riemann surface of genus $g\geq 0$, where $s_1,\, ...\, ,s_l$ (with $l\geq 0$) are $l$ removed points, and where $U$ is an arbitrary compact connected Lie group. We then give an interpretation of this Lagrangian submanifold in terms of decomposable representations. Finally, we study more carefully the relation between decomposable representations and irreducible representations.

\subsection{A Lagrangian submanifold of the moduli space $\Mgl$}

Recall from definition \ref{def_strat_lag_sub} that a \emph{stratified Lagrangian subspace} of a stratified symplectic space is a non-empty subset whose intersection with a given stratum, if non-empty, is a Lagrangian submanifold of the considered stratum. In section \ref{decomp}, we introduced, assuming the presence of an involutive automorphism of maximal rank $\tau$ of $U$ (see definition \ref{max_rk_inv}), a notion of \emph{decomposable representation} of $\pi_{g,l}$ into $U$ (see definition \ref{decomp_rep}) and then showed in section \ref{charac} that there exists an involution $\beta$ (see definition \ref{def_beta}) on the quasi-Hamiltonian space $$M=\Mtot$$ (where $\calC_1,\, ...\, ,\calC_l$ are $l$ conjugacy classes of $U$) such that a representation $$(a,b,c):=\abc$$ of $\pi_{g,l}$ into $U$ is decomposable if and only if $\beta(a,b,c)=\pphi.(a,b,c) $ for some $\pphi\in Fix(\taum)$, where $\taum$ is the involution $\taum(u):=\tau(u^{-1})$ of $U$ (see theorem \ref{charac_thm}). Recall that the representation space we are interested in is the quasi-Hamiltonian quotient $$\Mgl=M//U=\qhamquot$$ where $\mu$ is the momentum map
\begin{eqnarray*}
\mu:\Mext & \longto & U\\
\abc & \longmapsto & \relabc
\end{eqnarray*}
In order to find a stratified Lagrangian subspace of $\Mgl$, we will now use the considerations of subsection \ref{anti-symp_inv}. More precisely, we will show that $\beta$ induces an involution $\bhat$ on $M//U$ and that $Fix(\bhat)$ is a stratified Lagrangian subspace of $\Mgl$. 
\begin{theorem}[A stratified Lagrangian subspace of $\Mgl$]\label{lag_nature_thm}
For any compact connected Lie group $U$, the involution $\beta: M \to M$ induces an involution $\bhat$ on the moduli space $\Mgl=\Mod$ of representations of the surface group $\pi_{g,l}=\pigl$ into $U$ whose fixed-point set $Fix(\bhat)$ is a stratified Lagrangian subspace of $\Mgl$. In particular, for a generic choice of conjugacy classes $(\calC_j)_{1\leq j\leq l}$, the set $Fix(\bhat)$ is a Lagrangian submanifold of the smooth symplectic manifold $\Mgl$ ($l\geq 1$).
\end{theorem}
\begin{proof}
Compatibility of our $\beta$ with the diagonal conjugacy action of $U$ on $M$, and with the momentum map $\mu$ above, was verified in proposition \ref{beta_anti-inv}. Consequently, by proposition \ref{lag_locus}, it induces an involution $\bhat$ on $\Mgl$. We also showed in proposition \ref{beta_anti-inv} that $\beta$ reverses the $2$-form defining the quasi-Hamiltonian structure on $M=\Mtot$ (this was a consequence of the procedure for constructing form-reversing involutions on product spaces given in proposition \ref{anti-inv}) and that $\mu(Fix(\beta))$ intersects the connected component of $Fix(\taum)$ containing $1\in U$. Since the center of $U$ acts trivially on $M=\Mtot$ (the action of $U$ being the diagonal conjugacy action), theorem \ref{inv_arbitrary} applies and shows that we have $Fix(\beta)\cap\mu^{-1}(\{1\})\not=\emptyset$, and consequently $Fix(\bhat)\not=\emptyset$. Proposition \ref{fixbetadescription} then shows that if $Fix(\bhat)$ intersects a given stratum $X_j=(\fiber\cap M_{K_j})/L_{K_j}$ of $\Mgl=M//U$, then this intersection is the fixed-point set of an involution of the symplectic manifold $X_j$: $Fix(\bhat)\cap X_j = Fix(\bhat|_{X_j})$. Since $\beta^*w=-\w$ on $M$, $\bhat|_{X_j}$ is anti-symplectic. Therefore, if its fixed-point set is non-empty, it is a Lagrangian submanifold of $X_j$.
\end{proof}

\subsection{An interpretation of $Fix(\bhat)$ in terms of decomposable representations}\label{interpretation}

\newcommand{\calE}{\mathcal{E}}
Since the involution $\beta$ on $M=\Mtot$ providing the above example of a stratified Lagrangian subspace in $\Mgl=M//U$ characterizes decomposable representations of $\pi_{g,l}$ into $U$ (see theorem \ref{charac_thm}), it is natural to try to understand the relationship between the stratified Lagrangian subspace $Fix(\bhat)$ and the set of equivalence classes of decomposable representations of $\pi_{g,l}$. By proposition \ref{equiv_decomp_rep}, a representation $(a,b,c)$ is decomposable if and only if $u.(a,b,c)$ is decomposable for any $u\in U$. We may therefore denote by $\calE\subset\Mgl$ the set of equivalence classes of decomposable representations of $\pi_{g,l}$. By theorem \ref{charac_thm}, we have: $\calE\subset Fix(\bhat)$. We can in fact be more precise. Recall from proposition \ref{fixbetadescription} that we have the following description of $Fix(\bhat)$: $$Fix(\bhat)=\bigsqcup_{j\in J\ |\ \tau(K_j)\sim K_j} Fix(\bhat|_{X_j})$$ which we may rewrite under the following form: $$Fix(\bhat)=\Big(\bigsqcup_{j\in J_1} Fix(\bhat|_{X_j})\Big) \quad\bigsqcup\quad \Big(\bigsqcup_{j\in J_2} Fix(\bhat|_{X_j})\Big)$$ where $J_1=\{j\in J\ |\ \tau(K_j)=w_jK_jw_j^{-1}\ \mathrm{for\ some}\ w_j\in Fix(\taum)\}$ and $J_2=\{j\in J\ |\ \tau(K_j)\sim K_j\} \backslash J_1$. Then, for each $j\in J_1$, there exists, by proposition \ref{symmetric_case},  an involution $\beta_{K_j}:=w_j^{-1}.\beta:M_{K_j}\to M_{K_j}$ and an involutive automorphism $\tau_{K_j}:=Ad\ w_j^{-1}\circ \tau:U\to U$ such that $\beta_{K_j}$ is compatible with $\tau_{K_j}$ and the momentum map $\mu_j: M_{K_j}\to U$, and such that $Fix(\bhat|_{X_j})=Fix(\widehat{\beta_{K_j}})$. We then have, still by theorem \ref{charac_thm}: $\calE\subset \sqcup_{j\in J_1}Fix(\bhat|_{X_j})$. Observe here that each $Fix(\bhat|_{X_j})$ is a union of connected components of $Fix(\bhat)$ and that it is a  (smooth) Lagrangian submanifold of $X_j=(\fiber\cap M_{K_j})/L_{K_j}$. The natural question to ask is then: does one have the converse inclusion $\sqcup_{j\in J_1}Fix(\bhat|_{X_j})\subset \calE$? We now point out a number of cases in which this is true.
\begin{proposition}\label{distinct_eigenvalues}
Take $g\geq 0$ and $l\geq 1$. Let $U=U(n)$ and $M=\Mtot$ with $\calC_l\subset U(n)$ a conjugacy class defined by pairwise distinct eigenvalues. Then:
\begin{enumerate}
\item[(i)] for any isotropy group $K=U_{(a,b,c)}$ of a point $(a,b,c)\in M$, the group $\tau(K)$ is conjugate to $K$ by a symmetric element of $U$ (that is, $J_1=J$).
\item[(ii)] if $\beta(a,b,c)=\pphi.(a,b,c)$ for some $\pphi \in U$ then necessarily $\pphi\in Fix(\taum)$.
\item[(iii)] The set $\calE$ of equivalence classes of decomposable representations of $\pi_{g,l}$ into $U$ is exactly $Fix(\bhat)$.
\end{enumerate}
\end{proposition}
\begin{proof}\newcommand{\calZ}{\mathcal{Z}}
The proof will show that the above proposition is true for an arbitrary compact connected Lie group $U$ with $\calC_l$ the conjugacy class of an element $d_l\in U$ whose centralizer $\calZ(d_l)$ is a \emph{maximal torus} $T$ of $U$. Since $\tau$ is assumed to be of maximal rank, we may assume that this maximal torus is contained in $Fix(\taum)$.
\begin{enumerate}
\item[(i)] Let $K$ be the isotropy group of some $(a,b,c)\in M$. By assumption on $\calC_l$, there exists a $\psi\in U$ such that $d_l:=\psi c_l\psi^{-1}\in T$ and $\calZ(d_l)=T$ (that is, $d_l$ is diagonal with pairwise distinct eigenvalues). Take then $k\in K$. In particular $kc_l k^{-1}=c_l$ hence $\psi k\psi^{-1}\in \calZ(\psi c_l\psi^{-1}) =\calZ(d_l)=T$. Since $T\subset Fix(\taum)$, this implies $\taum(\psi k\psi^{-1})=\psi k\psi^{-1}$, hence $\taum(k)=(\taum(\psi)\psi)k(\taum(\psi)\psi)^{-1}$ and therefore $\taum(K)=(\taum(\psi)\psi)K(\taum(\psi)\psi)^{-1}$, with $(\taum(\psi)\psi)\in Fix(\taum)$.
\item[(ii)] Suppose that $\beta(a,b,c)=\pphi.(a,b,c)$ for some $\pphi \in U$ (that is, $[a,b,c]\in Fix(\bhat)$). Then in particular $\taum(c_l)=\pphi c_l\pphi^{-1}$. Hence $\taum(\psi^{-1}d_l\psi)=\pphi(\psi^{-1}d_l\psi)\pphi^{-1}$, where $\psi$ and $d_l$ are defined as in (i). Then $\taum(d_l)=(\tau(\psi)\pphi\psi^{-1})d_l(\tau(\psi)\pphi\psi^{-1})^{-1}$. But $\taum(d_l)=d_l$, hence $$(\tau(\psi)\pphi\psi^{-1})\in \calZ(d_l)=T$$ and therefore $\pphi\in \taum(\psi)T\psi\subset Fix(\taum)$ since $T\subset Fix(\taum)$.
\item[(iii)] We already know that $\calE\subset Fix(\bhat)$. The above shows that if $[a,b,c]\in Fix(\bhat)$ then $\beta(a,b,c)=\pphi.(a,b,c)$ with $\pphi\in Fix(\taum)$. By theorem \ref{charac_thm}, $(a,b,c)$ is then decomposable and we therefore have $Fix(\bhat)\subset \calE$.
\end{enumerate}
\end{proof}
In particular, one gets the following result as a corollary of the above:
\begin{corollary}\label{generic_case_result}
For a generic choice of conjugacy classes $\calC_1, ..., \calC_l$, the moduli space $\Mgl$ ($l\geq 1$) is a smooth symplectic manifold and the set $\calE$ of equivalence classes of decomposable representations is a smooth Lagrangian submanifold of $\Mgl$ which coincides with $Fix(\bhat)$.
\end{corollary}
Likewise, there is another special case in which $\calE=Fix(\bhat)$:
\begin{lemma}[\cite{Foth_loci}]
If $U$ acts freely on $\mu^{-1}(\{1\})$ then $\beta(a,b,c)=\pphi.(a,b,c)$ for some $\pphi\in U$ implies $\pphi\in Fix(\taum)$.
\end{lemma}
\begin{proof}
This is a general lemma on Hamiltonian and quasi-Hamiltonian quotients. One has $\beta(a,b,c)=\pphi.(a,b,c)$ hence $$(a,b,c)=\beta^2(a,b,c)=\beta(\pphi.(a,b,c))=\tau(\pphi).\beta(a,b,c)=\tau(\pphi)\pphi.(a,b,c).$$ Since the action of $U$ is free, this implies $\tau(\pphi)\pphi=1$ hence $\taum(\pphi)=\pphi$.
\end{proof}
Then, as in point (iii) of proposition \ref{distinct_eigenvalues}, the condition $[a,b,c]\in Fix(\bhat)$ implies $\beta(a,b,c)=\pphi.(a,b,c)$ with $\pphi\in Fix(\taum)$, hence by theorem \ref{charac_thm} $[a,b,c]\in\calE$. As an example of a group $U$ acting freely on $\mu^{-1}(\{1\})$, consider $U=PU(n)$ acting on $\Mgl$ with $l\geq 1$ and for a generic choice of conjugacy classes: in this case, one has $U_{(a,b,c)}=\mathcal{Z}(PU(n))=\{1\}$ for any representation $(a,b,c)\in\mu^{-1}(\{1\})$. In general though, the action of $U$ on $\fiber$ is not free for representations spaces of surface groups. But since the strata of $M//U$ are of the form $X_j=(\fiber\cap M_{K_j})/L_{K_j}$ with $L_{K_j}$ acting freely on $M_{K_j}$, one gets similarly:
\begin{proposition}\label{almost_there}
Take $j\in J_1$ and $[a,b,c]\in Fix(\bhat|_{X_j})=Fix(\widehat{\beta_{K_j}})\subset X_j$. Then $\beta_{K_j}(a,b,c)=\pphi_j.(a,b,c)$ for some $\pphi_j\in \calN(K_j)$ satisfying $\taum(\pphi_j)=\pphi_j k_j$, $k_j\in K_j$.
\end{proposition}
\begin{proof}
Assume that the equivalence class $[a,b,c]$ of a representation $(a,b,c)\in\mu^{-1}(\{1\})$ lies in $Fix(\bhat|_{X_j})$ for some $j\in J_1$. In particular, this means that we can assume that  the isotropy group of $(a,b,c)$ is exactly $K_j$: $(a,b,c)\in M_{K_j}$. Since $j\in J_1$, $\tau(K_j)$ is conjugate to $K_j$ by a symmetric element $w_j\in Fix(\taum)$ and we recall from proposition \ref{symmetric_case} that $\beta_{K_j}=w_{K_j}^{-1}.\beta:M_{K_j}\to M_{K_j}$ and $\tau_{K_j}=Ad\ w_j^{-1}\circ \tau$. Then, the fact that $[a,b,c]\in Fix(\bhat|_{X_j})=Fix(\widehat{\beta_{K_j}})$ yields $\beta_{K_j}(a,b,c)=[\pphi_j].(a,b,c)$ for some $[\pphi_j]\in L_{K_j}=\calN(K_j)/K_j$.We also denote by $\tau_{K_j}$ the involution induced by $\tau_{K_j}$ on $L_{K_j}$ (see lemma \ref{NK_stable}). We then have $(a,b,c)=\beta_{K_j}^2(a,b,c)=\beta_{K_j}([\pphi_j].(a,b,c))=\tau_{K_j}([\pphi_j]).\beta_{K_j}(a,b,c)=\tau_{K_j}([\pphi_j])[\pphi_j].(a,b,c)$. Since the action of $L_{K_j}$ on $M_{K_j}$ is free, one then has $\tau_{K_j}^{-}([\pphi_j])=[\pphi_j]$, which means that $\taum_{K_j}(\pphi_j)=\pphi_j k_j$ for some $k_j\in K_j$. Finally, since $K_j$ stabilizes $(a,b,c)$, one has $\beta_{K_j}(a,b,c)=\pphi_j.(a,b,c)$.
\end{proof}
\begin{remark}
If one had $k_j=1$ in the above, one would get $\pphi_j\in Fix(\taum_{K_j})$ which
by definition of $\beta_{K_j}$ and $\tau_{K_j}$ would translate, using the fact that $w_j\in Fix(\taum)$, to $\beta(a,b,c)=(w_j\pphi_j).(a,b,c)$ and $\taum(w_j\pphi_j)=w_j\pphi_j$. Theorem \ref{charac_thm} would then show that $(a,b,c)$ is decomposable.
\end{remark}
\begin{remark}
Proposition \ref{almost_there} combined with theorem \ref{charac_thm} also says that if one thinks of $(a,b,c)\in M_{K_j}$ (which implies $a_i,b_i,c_j \in \calN(K_j)$) as a representation of $\pi_{g,l}$ into $L_{K_j}$ then it is a decomposable one. Thus, in some sense, the strata of the stratified Lagrangian subspace $Fix(\bhat)\subset \Mgl$ consist of (equivalence classes of) decomposable representations.
\end{remark}\newcommand{\calZ}{\mathcal{Z}}

\subsection{Decomposable representations and irreducible representations}\label{decomp_and_irred}

We end our study of decomposable representations with remarks on the relation this notion bears to the notion of irreducible representation. If $U$ is an arbitrary compact connected Lie group and $\calZ(U)$ is the center of $U$, the elements of $\fiber \cap M_{\calZ(U)}$ are called \emph{irreducible representations} of $\pi_{g,l}$ into $U$ (if $U\subset Gl(V)$ is a group of linear transformations, those are indeed irreducible representations in the usual sense). Since $\calZ(U)$ is the smallest possible isotropy group of a point $(a,b,c)$ in $M=\Mtot$, equivalence classes of irreducible representations constitute the stratum of maximal dimension of $\Mgl=M//U$ (see \cite{Gol}). This stratum is called the \emph{principal stratum}. We then observe:
\begin{proposition}
If $(a,b,c)\in\fiber\cap M_{\calZ(U)}$ is an irreducible representation of $\pi_{g,l}$, then so is $\beta(a,b,c)$. More generally: $\beta(M_{\calZ(U)})=M_{\calZ(U)}$.
\end{proposition}
\begin{proof}
By proposition \ref{stab_beta_x}, one has $U_{\beta(a,b,c)}=\tau(U_{(a,b,c)})=\tau(\calZ(U))$. But then $\tau(\calZ(U))=\calZ(U)$ since $\tau$ is an involutive automorphism of $U$. The proposition follows.
\end{proof} 
The natural question to ask is then: does the principal stratum of $M//U$ always contain decomposable representations? In other words, do there always exist irreducible representations which are decomposable? In \cite{FW}, Falbel and Wentworth showed that, when $g=0$, $l\geq1$ and $U=U(n)$, there always exist irreducible decomposable representations. Their proof relies on an explicit computation of dimension in that particular case. Here, we treat another special case.
\begin{proposition}\label{ppal_stratum}
If $g\geq 1$ and $l=0$, then there always exist decomposable representations in the principal stratum of $\Mgl=\mathcal{M}_{g,0}$.
\end{proposition}
\begin{proof}
To show the existence of decomposable representations in the principal stratum of $\Mgl$, it is sufficient to apply theorem \ref{existence_thm_III} to the quasi-Hamiltonian $U$-space $M_{\calZ(U)}$. And to do that, all one has to prove is that there exists a fixed point of $\beta|_{M_{\calZ(U)}}$ whose image under $\mu$ lies in the connected component of $1$ in $Fix(\taum)$. Recall that here $M=(U\times U)^g$. Take $(a_1,b_1,\, ...\, , a_g,b_g)=(r_1,\tau(r_1),\, ...\, ,r_g,\tau(r_g))$ where $r_i\in \calZ(U)$ for $i\geq 2$ and $r_1\in U$ such that $U_{r_1}=\calZ(U)$. Then $(a,b)\in M_{\calZ(U)}$ since $U_{r_1}=\calZ(U)$ and $\beta(a,b)=(a,b)$ since the $r_i$ and the $\tau(r_i)$ commute to everything for $i\geq 2$. Moreover, $\mu(a,b)=r_1\tau(r_1)r_1^{-1}\tau(r_1)^{-1}=(r_1\tau(r_1))\taum(r_1\tau(r_1))$ does indeed lie in the connected component of $1$ in $Fix(\taum)$.
\end{proof}
Thus, if either ($g=0$ and $l\geq 1$) or ($g\geq 1$ and $l=0$), then there exist decomposable representations in the principal stratum of $\Mgl$, that is, decomposable representations of $\pi_{g,l}$ into $U$ which are irreducible. We do not know whether this is true for abitrary $g$ and $l$.  However, when ($g\geq 1$ and $l\geq 1$) and for a generic choice of conjugacy classes $(\calC_j)_{1\leq j\leq l}$, the moduli space $\Mgl$ ($l\geq 1$) is a smooth manifold, thus equal to the above-mentioned principal stratum: if $M_{\mathcal{Z}(U)}\not=\emptyset$ then all representations are irreducible. Since by theorem \ref{existence_thm_III} there always exist decomposable representations, these are irreducible. By corollary \ref{generic_case_result}, the set of equivalence classes of decomposable representations is then exactly $Fix(\bhat)$, a Lagrangian submanifold of $\Mgl$.

\section{Application: a compact analogue of the Thompson problem}\label{application}

In this section, we give an application of our results to a matrix problem: given $l$ conjugacy classes $\calC_1,\, ...\, ,\calC_l$ of the unitary group $U(n)$, when do there exist $l$ unitary matrices $A_1,\, ...\, ,A_l\in U(n)$ such that $A_j^tA_j\in \calC_j$ (where $A_j^t$ designates the transpose matrix of $A_j$) and $A_1...A_l=1$?\\
First, let us explain why we call this problem a compact analogue of Thompson's problem. Starting from the compact connected Lie group $U=U(n)$, the complexified group is $G:=U^{\C}=Gl(n,\C)$. The involutive automorphism $\theta: g\in Gl(n,\C)\mapsto (\bar{g}^t)^{-1}$ satisfies $Fix(\theta)=U(n)$. The \emph{non-compact dual} of $U(n)$ in $Gl(n,\C)$ is the group $H=Gl(n,\R)$: it satisfies $G=Fix(\tau)$ with $\tau:g\in Gl(n,\C)\mapsto \bar{g}$ and $\tau\circ\theta=\theta\circ\tau$. Since $\tau\circ\theta=\theta\circ\tau$, the involutive automorphism $\theta$ restricts to $Gl(n,\R)$, and the involution $\theta^{-}(h):=\theta(h^{-1})$ is equal to $\theta^{-}:h\in Gl(n,\R) \mapsto h^t$. We can then ask the following question: given $l$ conjugacy classes $\calC_1,\, ...\, ,\calC_l$ of $Gl(n,\R)$, when does there exists $l$ invertible real matrices $h_1,\, ...\, ,h_l\in Gl(n,\R)$ satisfying $\theta^{-}(h_j)h_j=h_j^th_j\in \calC_j$ and $h_1...h_l=1$? This problem is known as the real Thompson problem: the condition $h_j^th_j\in\calC_j$ says that $h_j$ has prescribed \emph{singular values}. The complex Thompson problem is formulated in $Gl(n,\C)$, asking when there exist $g_1,\, ...\, ,g_l\in Gl(n,\C)$ satisfying $\theta^{-}(g_j)g_j=g_j^*g_j\in \calC_j$ and $g_1...g_l=1$. The \emph{compact analogue} of this problem is the corresponding question when one replaces the involution $\theta$ on $H$ with the involution $\tau$ on the \emph{compact group} $U$. We now recall the following theorem of Klyachko in \cite{Klyachko}:
\begin{theorem}[\cite{Klyachko}]\label{Klyachko-AMW}
Consider $\lambda_1,\, ...\, ,\lambda_l\in\R^n$. Then the following
statements are equivalent:
\begin{enumerate}
\item[(i)] There exist $l$ invertible complex
matrices $g_1,\, ...\, ,g_l\in Gl(n,\C)$ such that:
$$\mathrm{Spec}\,(g_j^*g_j)=\exp(\lambda_j)\quad and \quad g_1...g_l=1.$$
\item[(ii)] There exist $l$ complex Hermitian matrices $H_1,\, ...\,
,H_l \in \mathcal{H}(n)$ such that: $$\mathrm{Spec}\, H_j =\lambda_j\quad
and\quad H_1+\cdots + H_l=0.$$
\item[(iii)] There exist $l$ invertible real
matrices $h_1,\, ...\, ,h_l\in Gl(n,\R)$ such that:
$$\mathrm{Spec}\,(h_j^th_j)=\exp(\lambda_j)\quad and\quad h_1...h_l=1.$$
\item[(iv)] There exist $l$ real symmetric matrices $S_1,\, ...\,
,S_l \in \mathcal{S}(n)$ such that: $$\mathrm{Spec}\, S_j =\lambda_j\quad
and\quad S_1+\cdots + S_l=0.$$
\end{enumerate}
\end{theorem}
The equivalence of (i) and (ii) is Thompson's conjecture. It was first proved by Klyachko in \cite{Klyachko} and then proved using symplectic geometry in \cite{AMW}. Statements (iii) and (iv) are real versions of (i) and (ii), respectively. We refer to \cite{AMW} to see how the equivalence of (ii) of (iv) relies on a real convexity theorem for Lie-algebra-valued momentum maps. For a thorough treatment of the Lie-theoretic approach to Thompson's conjecture, we refer to the work of Evens and Lu in \cite{Evens-Lu}. Finally, we refer to \cite{Klyachko-Herm} for necessary and sufficient conditions on the $\lambda_j$ for problem (ii) (and therefore any of the other three) to have a solution. These conditions are linear inequalities satisfied by the $\lambda_j$.\\
For the compact Thompson problem that we formulated at the beginning of this section, the complex version should be formulated in the compact group $U\times U\subset G\times G=G^{\C}$. For further details in this direction, we refer to \cite{Sch_CJM}, where we also explain formally what the expression \emph{real version} of a problem means. We will now prove a theorem that provides an answer to the question asked at the beginning of this section:
\begin{theorem}[An application to a matrix problem]\label{compact_Thompson}
Consider $\lambda_1, ...\, ,$ $\lambda_l\in\R^n$. Then the following
statements are equivalent:
\begin{enumerate}
\item[(i)] There exist $l$ unitary
matrices $u_,\, ...\, ,u_l\in U(n)$ such that:
$$\mathrm{Spec}\,u_j=\exp(i\lambda_j)\quad and \quad u_ l...u_l=1.$$
\item[(ii)] There exist $l$ unitary matrices $A_1,\, ...\, ,A_l \in U(n)$ such that: $$\mathrm{Spec}\, (A_j^tA_j) =\exp(i\lambda_j)\quad
and\quad A_1...A_l=1.$$
\end{enumerate}
\end{theorem}
\begin{proof}
Let $\calC_j$ be the conjugacy class of $\exp(i\lambda_j)\in U(n)$. Assume first that $A_1.\, ...\, ,A_l\in U(n)$ satisfy $A_j^tA_j\in \calC_j$ and $A_1...A_l=1$. Set $u_l=A_l^tA_l$, $u_{l-1}=A_l^t(A_{l-1}^tA_{l-1})(A_l^t)^{-1}$, ...\, , and $u_1=(A_2...A_l)^t(A_1^tA_1)((A_2...A_l)^t)^{-1}$. Then $u_j\in\calC_j$ for all $j$ and
\begin{eqnarray*}
u_1...u_l & = & (A_2...A_l)^t(A_1)^tA_1(A_2...A_l)\\
& = & (A_1...A_l)^t(A_1...A_l)\\
& = & 1
\end{eqnarray*}
which proves that (ii) implies (i).\\
Conversely, consider $u_1,\, ...\, ,u_l\in U(n)$ satisfying $u_j\in\calC_j$ and $u_1...u_l=1$. This means that the momentum map
\begin{eqnarray*}
\mu:\pconj & \longto & U(n)\\
(u_1,\, ...\, ,u_l) & \longmapsto & u_1...u_l
\end{eqnarray*}
has a non-empty fiber above $1\in U(n)$. The Lie group $U(n)$ is endowed with the involution $\tau(u)=\overline{u}$ (therefore $\taum(u)=u^t$). Since the center of $U(n)$ acts trivially on the quasi-Hamiltonian space $\pconj$, theorem \ref{inv_arbitrary} applies and shows that $Fix(\beta)\cap\mu^{-1}(\{1\})\not=\emptyset$, where $\beta$ is the involution on $\pconj$ introduced in definition \ref{def_beta}. Consider then $(w_1,\, ...\, ,w_l)\in Fix(\beta)\cap\mu^{-1}(\{1\})$. Then $\beta(w_1,\, ...\, ,w_l)=(w_1,\, ...\, ,w_l)$, that is:
\begin{eqnarray*}
w_l^t & = & w_l\\
w_l^t w_{l-1}^t (w_l^t)^{-1} & = & w_{l-1}\\
& \vdots & \\
\big(w_2...w_l\big)^tw_1^t\big((w_2...w_l)^t\big)^{-1} & = & w_1
\end{eqnarray*}
Since $Fix(\taum)\subset U(n)$ is connected (every symmetric unitary matrix is of the form $w=\exp(iS)$ where $S$ is a real symmetric matrix), we can write $w_l=A_l^tA_l$ for some $A\in U(n)$ (take for instance $A=\exp(i\frac{S}{2}))$. Using the above equations, we can then write $A_l^tA_lw_{l-1}^tA_l^{-1}(A_l^t)^{-1}=w_{l-1}$, that is: $$\big((A_l^t)^{-1}w_{l-1}A_l^t\big)^t=(A_l^t)^{-1}w_{l-1}A_l^t$$ and we can therefore write $$(A_l^t)^{-1}w_{l-1}A_l^t=A_{l-1}^tA_{l-1}$$ for some $A_{l-1}\in U(n)$. Continuing like this, we obtain, for all $j\in\{2,\, ...\, ,n\}$:
\begin{eqnarray}
\big((A_{j+1}...A_l)^t\big)^{-1}w_j\big(A_{j+1}...A_l\big)^t & =  & A_j^tA_j.\label{symrel}
\end{eqnarray}
In particular, $A_j^tA_j\in\calC_j$ for all $j\geq 2$. We then set $A_1:=(A_2...A_l)^{-1}$. Then:
\begin{eqnarray*}
& & \big(A_2...A_l\big)^t\big(A_1^tA_1\big)\big((A_2...A_l)^t\big)^{-1}\\
& = & \big(A_2...A_l\big)^t\big((A_2...A_l)^{-1}\big)^t \big(A_2...A_l)^{-1} \big((A_2...A_l)^t\big)^{-1}\\
& = & \big((A_2...A_l)^t(A_2...A_l)\big)^{-1}\\
& = & \big(A_l^t...A_3^t(A_2^tA_2)A_3...A_l\big)^{-1}\\
\mathrm{\big(using\ (}\ref{symrel}\mathrm{)\big)} & = & \big(\underbrace{A_l^t...A_3^t\big((A_3...A_l)^t\big)^{-1}}_{=1}\, w_2 \, (A_3...A_l)^tA_3...A_l\big)^{-1}\\
\mathrm{\big(by\ induction\big)} & = & \big(w_2w_3...w_l)^{-1}\\
& = & w_1
\end{eqnarray*}
since $w_1...w_l=1$. In particular, $A_1^tA_1$ is conjugate to $w_1$ and therefore $A_1^tA_1\in\calC_1$. Since $A_1...A_l=1$ by definition of $A_1$, this shows that (i) implies (ii). 
\end{proof}
Observe that the most important part of the proof is the use of theorem \ref{inv_arbitrary} to ensure the existence of some $(w_1,\, ...\, ,w_l)\in Fix(\beta)\cap\mu^{-1}(\{1\})$ (which, in our terminology, means that there exists decomposable representations of the fundamental group $\piS$ of the punctured sphere). As we have seen in section \ref{existence}, this ultimately relies on the real convexity theorem for group-valued momentum maps \ref{convexity_thm}, just as the proof of theorem \ref{Klyachko-AMW} given in \cite{AMW} ultimately relies on a real convexity theorem for Lie-algebra-valued momentum maps. Finally, as for Thompson's problem, precise conditions on the $\lambda_j$ above for statement (i) to be true are already known (we refer for instance to \cite{AW}): these conditions are linear inequalities in the $\lambda_j$.

\begin{ack}
This paper was written during my post-docto\-ral stay at Keio University in Yokohama, Japan. The ideas developed here find their origin in \cite{Sch_CJM} (where the case of the punctured sphere group was treated). During the past three years, I have received encouragement, suggestions or been indicated a certain book or article by many people whom I would like to thank for their help: Alan Weinstein, Johannes Huebschmann, Sam Evens, Philip Boalch, Pierre Sleewaegen, Anton Alekseev, Eckhard Meinrenken. In addition to that, I have been very fortunate to have had regular interlocutors to discuss this problem with: Elisha Falbel, Jiang-Hua Lu, Richard Wentworth, Pierre Will. I would like to thank them sincerely for their interest and their patience. Finally, I would like to thank Yoshiaki Maeda, my host in Keio University, for giving me the opportunity to work in the excellent conditions provided by the Keio mathematics department.
\end{ack}

\end{document}